\documentclass[11pt]{article}
\usepackage{fullpage}
\usepackage{subfig}
\usepackage[usenames,dvipsnames]{xcolor}
\usepackage[colorlinks,citecolor=blue,linkcolor=BrickRed]{hyperref}
	\usepackage{latexsym,graphicx,epsfig,color}
\usepackage{amsfonts,amssymb,amsmath,amsthm,amstext}
\usepackage{libertine}
\usepackage[libertine]{newtxmath}
\usepackage{enumitem}
\setitemize{itemsep=2pt,topsep=0pt,parsep=0pt}
\usepackage{url,setspace}
\usepackage{multirow}
\usepackage{rotating}
\usepackage{makeidx}
\usepackage{tikz}
\usepackage{accents}
\usepackage{xspace}
\usepackage{algorithm,algpseudocode}
\usepackage{bm}
\usepackage{changepage} 	%\begin{adjustwidth}{0.5cm}{} \end{adjustwidth}
\usepackage{thmtools,thm-restate}
\usepackage{nicefrac}
\usepackage{cleveref}
\usepackage{import}
\usepackage[toc]{appendix}

\newcommand{\IGNORE}[1]{}
\allowdisplaybreaks

\usetikzlibrary{decorations.markings}
\usetikzlibrary{arrows}
\tikzstyle{block}=[draw opacity=0.7,line width=1.4cm]
\tikzstyle{graphnode}=[circle, draw, fill=black!20, inner sep=0pt, minimum width=6pt]
\tikzstyle{point}=[circle, draw, fill=black!30, inner sep=0pt, minimum width=1pt]
\tikzstyle{input}=[rectangle, draw, fill=black!75,inner sep=3pt, inner ysep=3pt, minimum width=4pt]
\tikzstyle{unmatched}=[graphnode,fill=black!0]
\tikzstyle{shaded}=[graphnode,fill=black!20]
\tikzstyle{matched}=[graphnode,fill=black!100]  	
\tikzstyle{matching} = [ultra thick]
\tikzset{
    >=stealth',
    pil/.style={
           ->,
           thick,
           shorten <=2pt,
           shorten >=2pt,}
}
\tikzset{->-/.style={decoration={
  markings,
  mark=at position .5 with {\arrow{>}}},postaction={decorate}}}
\makeindex


\newtheorem{theorem}{Theorem}[section]
\newtheorem{claim}[theorem]{Claim}
\newtheorem{proposition}[theorem]{Proposition}
\newtheorem{lemma}[theorem]{Lemma}

\newtheorem{fact}[theorem]{Fact}
\newtheorem*{remark}{Remark}
\theoremstyle{definition}

\newtheorem{defn}[theorem]{Definition}

\newtheorem{rem}[theorem]{Remark}

\usepackage{thmtools,thm-restate} % Used to restate lemmas without new introducing new numbers. 

\newcommand{\mul}{\mathrm{mul}}
\newcommand{\cnt}{\mathrm{cnt}}

\newcommand{\calC}{\mathcal{C}}

\newcommand{\calF}{\mathcal{F}}
\newcommand{\calK}{\mathcal{K}}

\newcommand{\calP}{\mathcal{P}}
\newcommand{\calT}{\mathcal{T}}
\newcommand{\calI}{\mathcal{I}}

\def \naturals {\mathbb{N}}

\title{Powers of Hamilton cycles of high discrepancy are unavoidable}
	\author{Domagoj Brada\v{c}\thanks{Department of Mathematics, ETH Z\"{u}rich, 8092 Z\"{u}rich, Switzerland. Email: \texttt{domagoj.bradac@math.ethz.ch}}
}

\date{ \today}

\begin{document}
\maketitle
\thispagestyle{empty}

\setlength{\abovedisplayskip}{2pt}
\setlength{\belowdisplayskip}{2pt}

\begin{abstract}
The P\'osa-Seymour conjecture asserts that every graph on $n$ vertices with minimum degree at least $(1 - 1/(r+1))n$ contains the $r^{th}$ power of a Hamilton cycle. Koml\'os, S\'ark\"ozy and Szemer\'edi famously proved the conjecture for large $n.$ The notion of discrepancy appears in many areas of mathematics, including graph theory. In this setting, a graph $G$ is given along with a $2$-coloring of its edges. One is then asked to find in $G$ a copy of a given subgraph with a large discrepancy, i.e., with significantly more than half of its edges in one color. For $r \geq 2,$ we determine the minimum degree threshold needed to find the $r^{th}$ power of a Hamilton cycle of large discrepancy, answering a question posed by Balogh, Csaba, Pluh\'ar and Treglown. Notably, for $r \geq 3,$ this threshold approximately matches the minimum degree requirement of the P\'osa-Seymour conjecture.
\end{abstract}
\thispagestyle{empty}

\setcounter{page}{1}

% !TeX root = main.tex
% !TEX root = main.tex

\section{Introduction}
%\subsection{Graph Discrepancy Theory}

Classical discrepancy theory studies problems of the following kind: given a family of subsets of a universal set $\mathcal{U},$ is it possible to partition the elements of $\mathcal{U}$ into two parts such that each set in the family has roughly the same number of elements in each part? One of the first significant results in this area is a criterion for a sequence to be uniformly distributed in the unit interval proved by Hermann Weyl \cite{weyl1916gleichverteilung}. Since then, discrepancy theory has had wide applicability in many fields such as ergodic theory, number theory, statistics, geometry, computer science, etc. For a comprehensive overview of the field, see the books by Beck and Chen \cite{beck_chen_1987}, Chazelle \cite{chazelle2001discrepancy} and Matou\v{s}ek \cite{matousek1999geometric}.

This paper studies a problem in the discrepancy theory of graphs. To discuss the topic, we start with a definition.
\begin{defn}
  Let $G = (V, E)$ be a graph and $f \colon E \rightarrow \{-1, 1\}$ a labelling of its edges. Given a subgraph $H$ of $G,$ we define its discrepancy $f(H)$ as 
  \[ f(H) = \sum_{e \in E(H)} f(e). \]
  Furthermore, we refer to the value $|f(H)|$ as the absolute discrepancy of $H.$
\end{defn}

One of the central questions in graph discrepancy theory is the following. Suppose we are given a graph $G$ and a spanning subgraph $H.$ Does $G,$ for every edge labelling $f\colon E(G) \rightarrow \{-1, 1\},$ contain an isomorphic copy of $H$ of high absolute discrepancy with respect to $f?$ Erd\H{o}s, F\"uredi, Loebl and S\'{o}s \cite{erdHos1995discrepancy} proved the first result of this kind. They show that for large enough $n,$ given a tree on $n$ vertices $T_n$ with maximum degree $\Delta$ and a $\{-1, 1\}$-coloring of the edges of the complete graph $K_n,$ one can find a copy of $T_n$ with absolute discrepancy at least $c(n-1-\Delta),$ for some absolute constant $c > 0.$ 

A commonly studied topic in extremal combinatorics are \textit{Dirac-type problems} where one is given a graph $G$ on $n$ vertices with minimum degree at least $\alpha n$ and wants to prove that $G$ contains a copy of a specific spanning subgraph $H$. In the discrepancy setting it is natural to ask whether we can also find a copy of $H$ with large absolute discrepancy. Balogh, Csaba, Jing and Pluh\'ar studied this problem for spanning trees, paths and Hamilton cycles. Among other results, they determine the minimum degree threshold needed to force a Hamilton cycle of high discrepancy.
\begin{theorem}[Balogh, Csaba, Jing and Pluh\'ar \cite{balogh2020discrepancies}]
  Let $0 < c < 1/4$ and $n \in \naturals$ be sufficiently large. If $G$ is an $n$-vertex graph with 
  \[\delta(G) \geq (3/4 + c)n \]
  and $f\colon E(G) \rightarrow \{-1, 1\},$ then $G$ contains a Hamilton cycle with absolute discrepancy at least $cn/32$ with respect to $f.$ Moreover, if $4$ divides $n,$ there is an $n$-vertex graph with $\delta(G) = 3n/4$ and an edge labelling $f\colon E(G) \rightarrow \{-1, 1\}$ such that any Hamilton cycle in $G$ has discrepancy $0$ with respect to $f.$
\end{theorem}

Very recently, Freschi, Hyde, Lada and Treglown \cite{freschi2020note}, and independently, Gishboliner, Krivelevich and Michaeli \cite{gishboliner2020discrepancies} generalized this result to edge-colorings with more than two colors. 

A fundamental result in extremal graph theory is the Hajnal-Szemer\'edi theorem. It states that if $r$ divides $n$ and $G$ is a graph on $n$ vertices with $\delta(G) \geq (1 - 1/r)n,$ then $G$ contains a perfect $K_r$-tiling, i.e. its vertex set can be partitioned into disjoint cliques of size $r.$ Balogh, Csaba, Pluh\'ar and Treglown \cite{balogh2020discrepancy} proved a discrepancy version of this theorem.
\begin{theorem}[Balogh, Csaba, Pluh\'ar and Treglown \cite{balogh2020discrepancy}]
\label{disc_hajnal_thm}
Suppose $r \geq 3$ is an integer and  let $\eta >0$. Then there exists $n_0 \in \mathbb N$ and $\gamma >0$ such that the following holds.
Let $G$ be a graph on $n \geq n_0$ vertices where $r$ divides $n$ and where
$$\delta (G) \geq \left (1-\frac{1}{r+1}+\eta \right )n.$$
Given any function $f\colon E(G) \rightarrow \{-1,1\}$ there exists a perfect $K_r$-tiling $\mathcal T$ in $G$ so that
$$\Big| \sum _{e \in E(\mathcal T)} f(e) \Big| \geq \gamma n.$$
Moreover, if $2r(r+1)$ divides $n,$ there exists a graph $G$ on $n$ vertices with $\delta(G) = (1 - 1/(r+1))n$ and a $2$-coloring of its edges such that the discrepancy of any perfect $K_r$-tiling of $G$ is $0.$
\end{theorem}

The $r^{th}$ power of a graph $G$ is the graph on the same vertex set in which two vertices are joined by an edge if and only if their distance in $G$ is at most $r.$ The P\'osa-Seymour conjecture asserts that any graph on $n$ vertices with minimum degree at least $(1 - 1/(r+1))n$ contains the $r^{th}$ power of a Hamilton cycle. Koml\'os, S\'ark\"ozy and Szemer\'edi \cite{komlos1998proof} proved the conjecture for large $n.$ In \cite{balogh2020discrepancy} the authors posed the question of determining the minimum degree needed to force the $r^{th}$ power of a Hamilton cycle with absolute discrepancy linear in $n.$ Because the $r^{th}$ power of a Hamilton cycle contains a (almost) perfect $K_{r+1}$-tiling, they suggested the minimum degree required should be $(1 - 1/(r+2) + \eta)n,$ based on their result for $K_r$-tilings. We prove this value is correct for $r=2.$ However, we show that for $r \geq 3,$ a minimum degree of $(1 - 1/(r+1) + \eta)n,$ for arbitrarily small $\eta > 0,$ is sufficient, approximately matching the minimum degree required for finding any $r^{th}$ power of a Hamilton cycle. As far as the author knows, this is the first Dirac-type discrepancy result in which the threshold for finding a spanning subgraph of large discrepancy is the same, up to an arbitrarily small linear term, as the minimum degree required for finding any copy of the subgraph.

\begin{theorem}
  \label{main_thm}
  For any integer $r \geq 3$ and $\eta > 0,$ there exist $n_0 \in \naturals$ and $\gamma > 0$ such that the following holds. Suppose a graph $G$ on $n \geq n_0$ vertices with minimum degree $\delta(G) \geq (1 - 1/(r+1) + \eta)n$ and an edge coloring $f \colon E(G) \rightarrow \{-1, 1\}$ are given. Then in $G$ there exists the $r^{th}$ power of a Hamilton cycle $H^r$ satisfying
  \[ \Big| \sum_{e \in E(H^r)} f(e) \Big| \geq \gamma n. \]
\end{theorem}

Interestingly, the minimum degree needed for finding the $r^{th}$ power of a Hamilton cycle of large discrepancy is the same for $r \in \{1, 2, 3\}$ and equals $(\frac{3}{4} + \eta)n.$ The cases $r=1,3$ being resolved in \cite{balogh2020discrepancies} and by the previous theorem, respectively, we also prove this for $r=2.$

\begin{theorem}
  \label{r2_thm}
  For any $\eta > 0,$ there exist $n_0 \in \naturals$ and $\gamma > 0$ such that the following holds. Suppose a graph $G$ on $n \geq n_0$ vertices with minimum degree $\delta(G) \geq (3/4 + \eta)n$ and an edge coloring $f \colon E(G) \rightarrow \{-1, 1\}$ are given. Then in $G$ there exists the square of a Hamilton cycle $H^2$ satisfying
  \[ \Big|\sum_{e \in E(H^2)} f(e) \Big| \geq \gamma n. \]
\end{theorem}

These results are tight in the following sense. If we weaken the minimum degree requirement by replacing the term $\eta n$ with a sublinear term, then there are examples in which any $r^{th}$ power of a Hamilton cycle has absolute discrepancy $o(n).$

The paper is organised as follows. In Section 2 we introduce some notation and definitions and state previous results used in our proofs. Then we present lower bounds showing the tightness of our results in Section 3.  We give a short outline of the proofs in Section 4. The proofs are then divided into two sections. In Section 5 we adapt the proof by Koml{\'o}s, S{\'a}rk{\"o}zy and Szemer{\'e}di of an approximate version of the P\'osa-Seymour conjecture to our setting, while the rest of the argument is presented in Section 6.

% !TeX root = main.tex
% !TEX root = main.tex

\section{Preliminaries}

Most of the graph theory notation we use is standard in the literature and can be found in ~\cite{bondy2007graph}. Let $G$ be a graph. We use $N_G(v)$ to denote the neighbourhood of $v$ in $G.$ For a set of vertices $S,$ we use $G[S]$ to denote the subgraph of $G$ induced by $S.$ We use $N_G(S)$ for the common neighbourhood of $S,$ formally $N_G(S) = \{ x \in V(G) \, \vert \, xv \in E(G), \forall v \in S\}$ and we denote $\deg_G(S) = |N_G(S)|.$ Given an edge labelling $f \colon E(G) \rightarrow \{-1, 1\},$ we use $G_+$ to denote the graph containing all edges labelled $1$ and $G_-$ to denote the graph containing all edges labelled $-1.$ We write $N_+(v)$ for the set of $u$ in $N_G(v)$ such that $f(v, u) = 1$ and $N_-(v)$ for the set of $u$ in $N_G(v)$ such that $f(v, u) = -1.$ We further denote $\deg_+(v) = \deg_{G_+}(v)$ and $\deg_-(v) = \deg_{G_-}(v).$ For a vertex $v$ and a subset of vertices $U,$ we define $N(v, U) = N(v) \cap U$ and $\deg(v, U) = |N(v, U)|.$ 
We write $\stackrel{.}{\cup}$ for the union of disjoint sets. We use the terms edge labelling and edge coloring interchangeably. We sometimes omit the underlying graph when it is clear from the context.

Throughout the paper we allow cycles to have repeated vertices, unless explicitly stated they are simple. For a cycle $C,$ we denote by $\cnt_C(v)$ the number of occurences of $v$ in $C$ when viewed as a closed walk.
The $r^{th}$ power of a cycle $C$ is defined as the multigraph obtained by adding edges between every pair of vertices whose indices differ by at most $r$ in the cyclic order of $C.$ We use $\mul_H(e)$ for the multiplicity of an edge $e$ in a multigraph $H.$ Then, for a cycle $C = (v_1, v_2, \dots, v_m),$ its $r^{th}$ power $C^r$ is formally defined as the multigraph on the vertex set $\{v_1, v_2, \dots, v_m\}$ with the following edge multiplicities:
\[ \mul_{C^r}(xy) = \left| \big\{ (i, j) \; \vert \; i \in [m], j \in [r], \{ v_i, v_{i+j} \} = \{x, y\} \big\} \right|, \] where we denote $v_{m+i} = v_i$ for $1 \leq i \leq r.$
The $r^{th}$ power of a simple $(r+1)$-cycle will sometimes be referred to as an $(r+1)$-clique. Importantly, though, it has two copies of each edge.
Given an edge labelling $f \colon E(G) \rightarrow \{-1, 1\}$ and a cycle $C$ such that $C^r$ is a subgraph of $G,$ we define the discrepancy of $C^r$ in the natural way: $f(C^r) = \sum_{e \in E(C^r)} \mul_{C^r}(e) f(e).$ Here we slightly abuse notation in the following sense. We ignore edge multiplicities for the notion of graph containment (as the ambient graph is always simple). In other words, we only require the ambient graph to have one copy of each edge that has positive multiplicity in a given $r^{th}$ power of a cycle.

Similarly as in \cite{balogh2020discrepancy}, we define a $C^r$-template. Note that in the following definition we only allow short cycles.

\begin{defn}
  Let $F$ be a graph. A $C^r$-template of $F$ is a collection $\calF = \{C_1, C_2, \dots, C_s\}$ of not necessarily distinct cycles whose $r^{th}$ powers are subgraphs of $F.$ In a $C^r$-template each vertex appears the same number of times, that is, $\sum_{i=1}^s \cnt_{C_i}(v)$ is the same for all $v \in V(F).$ Moreover, we require that each cycle $C_i$ has length between $r+1$ and $10r^2.$ The discrepancy of a $C^r$-template is given as $f(\calF) = \sum_{i=1}^s f(C_i^r).$ 
\end{defn}

The notion of a $C^r$-tiling is obtained by adding the natural restriction that each vertex appears exactly once. 
\begin{defn}
  Let $F$ be a graph. A $C^r$-tiling $\calT$ of $F$ is a collection of simple cycles $\calT = \{C_1, C_2, \dots, C_s\}$ whose $r^{th}$ powers are subgraphs of $F$ and each vertex appears precisely once in these cycles. The length of each cycle is between $r+1$ and $10r^2.$ The discrepancy of a $C^r$-tiling is given as $f(\calT) = \sum_{i=1}^s f(C_i^r).$
\end{defn}
A $K_{r+1}$-tiling can be viewed as a $C^r$-tiling in which all cycles have length $r+1,$ that is, all tiles are $(r+1)$-cliques. We sometimes implicitly use this interpretation of a $K_{r+1}$-tiling.

We give names to special types of $k$-cliques with respect to an edge labelling $f.$

\begin{defn}
  We write $K_k^+$ for the $k$-clique with all edges labelled $1$ and $K_k^-$ for the $k$-clique with all edges labelled $-1.$ The $(K_k, +)$-star is the clique whose edges labelled $1$ induce a copy of $K_{1, k-1}.$ The root of this $K_{1, k-1}$ is called the \emph{head} of the $(K_k, +)$-star. We define the $(K_k, -)$-star and its head analogously.
\end{defn}

We write $\alpha \ll \beta \ll \gamma,$ if the constants can be chosen from right to left such that all calculations in our proof are valid. More precisely, $\alpha \ll \beta$ means there is a positive increasing function $f(\beta)$ such that for $\alpha \le f(\beta),$ all calculations in the proof are valid. This notion naturally extends to hierarchies of larger length as well. We omit floors and ceilings whenever they do not affect the argument.

In our proofs we use the famous Hajnal-Szemer\'edi theorem in the following form.
\begin{theorem}[Hajnal and Szemer\'edi \cite{hajnal_szemeredi}]
  \label{hajnal_szemeredi}
  Every graph $G$ whose order $n$ is divisible by $r$ and has minimum degree at least $(1 - 1/r)n$ contains a perfect $K_r$-tiling.
\end{theorem}

\subsection{The Regularity Lemma}
In the proof of the main result, we use a multicolored variant of Szemer\'edi's regularity lemma \cite{szemeredi1975regular}. Before stating the result, we define the relevant notions. The density of a bipartite graph $G$ with vertex classes $A$ and $B$ is defined as
\[ d(A, B) = \frac{e(A, B)}{|A||B|}. \]
Given $\varepsilon, d > 0,$ the graph $G$ is said to be $(\varepsilon, d)$-regular if $d(A, B) \geq d$ and for any $X \subseteq A, Y \subseteq B$ such that $|X| > \varepsilon |A|$ and $|Y| > \varepsilon |B|,$ we have $|d(A, B) - d(X, Y)| < \varepsilon.$ The graph $G$ is $(\varepsilon, \delta)$-super-regular if for every $X \subseteq A$ with $|X| > \varepsilon |A|$ and $Y \subseteq B$ with $|Y| > \varepsilon |B|,$ we have $d(X, Y) > \delta,$ and furthermore, $\deg(a) > \delta |B|$ for all $a \in A$ and $\deg(b) > \delta |A|$ for all $b \in B.$ Suppose we are given a graph $G$ with an edge labelling $f \colon E(G) \rightarrow \{-1, 1\}.$ Given disjoint $X, Y \subseteq V(G)$ we write $(X, Y)_{G_+}$ or $G_+[X, Y]$ for the bipartite graph with vertex classes $X, Y$ containing edges between $X$ and $Y$ labelled $1.$ Analogously, we define $(X, Y)_{G_-}$ and $G_-[X, Y]$ with respect to edges labelled $-1.$

We use a variant of the regularity lemma which is easily deduced from the multicolored version in \cite{komlos1996szemeredi}. 
\begin{lemma}
  \label{reg_lemma}
  For every $\varepsilon > 0$ and $\ell_0 \in \naturals$ there exists $L_0 =
  L_0(\varepsilon, \ell_0)$ such that the following holds. Let $d \in [0, 1)$
  and let $G$ be a graph on $n \geq L_0$ vertices with an edge coloring $f\colon
  E(G) \rightarrow \{-1, 1\}.$ Then, there exists a partition
  $(V_i)_{i=0}^\ell$, for some $\ell \in [\ell_0, L_0]$ of $V(G)$ and a
  spanning subgraph $G'$ of $G$ with the following properties:
  \begin{enumerate}[label=(\roman*), noitemsep, topsep=-5pt]
    \item $|V_0| \leq \varepsilon n$ and $|V_1| = |V_2| = \dots = |V_\ell|$,
    \item $\deg_{G'}(v) \geq \deg_G(v) - (2d + \varepsilon)n$ for every $v \in V(G);$ \label{reg_lemma_min_degree}
    \item $e(G'[V_i]) = 0$ for all $1 \leq i \leq \ell;$
    \item for all $1 \leq i < j \leq \ell$ and $\sigma \in \{+, -\}$, either $(V_i,
      V_j)_{G'}^\sigma$ is an $(\varepsilon, d)$-regular pair or $G_\sigma'[V_i, V_j]$ is
      empty;
  \end{enumerate}
\end{lemma}

The above partition $V_0, V_1, \dots, V_\ell$ will be called an \emph{$(\varepsilon, d)$-regular partition of $G$ with respect to $f$}. We call $V_1, \dots, V_\ell$ \textit{clusters} and $V_0$ the \textit{exceptional set}. We refer to $G'$ as the \textit{pure graph}. We define the \textit{reduced graph} $R$ of $G$ with parameters $\varepsilon, d, \ell_0 $ to be the graph whose vertices are $V_1, \dots, V_\ell$ and where $(V_i, V_j)$ is an edge if at least one of $(V_i, V_j)_{G'_+}$ and $(V_i, V_j)_{G'_-}$ is $(\varepsilon, d)$-regular. On the reduced graph $R,$ we define the edge coloring $f_R \colon E(R) \rightarrow \{-1, 1\}$ as follows:
\begin{equation}
\label{deffR}
f_R(V_i, V_j) =
\begin{cases}
  1, &\text{if } (V_i, V_j)_{G'_+} \text{ is } (\varepsilon, d)\text{-regular} \\
 -1, &\text{ otherwise.}
\end{cases}
\end{equation}
Note that if both $(V_i, V_j)_{G'_+}$ and  $(V_i, V_j)_{G'_-}$ are $(\varepsilon, d)$-regular, $f_R$ only records the former property.\\

\begin{rem} \label{rem_divisibility}
  When working with the reduced graph, it will sometimes be convenient for us to assume that the number of its vertices is divisible by $r+1$ in order to apply Theorem~\ref{hajnal_szemeredi}. However, this can easily be achieved by moving at most $r$ clusters to the exceptional set. Provided $r / \ell_0 \le \varepsilon$ (which we will always be able to guarantee), we move at most $\varepsilon n$ vertices to the exceptional set. It is easy to check that the new graph satisfies the properties given by Lemma~\ref{reg_lemma} with the same parameters except that $\varepsilon$ increases by at most a factor of $2.$ As we always choose $\varepsilon$ to be sufficiently small, this does not affect any of our arguments.
\end{rem}

We use the following simple fact about the reduced graph.
\begin{fact}
  \label{degRfact}
  Let $c > 0$ be a given constant and $G$ a graph on $n$ vertices such that $\delta(G) \geq cn.$ Let $R$ be the reduced graph obtained after applying Lemma \ref{reg_lemma} with parameters $\varepsilon, d$ and $\ell_0.$ Then $\delta(R) \geq (c - 2d - 2\varepsilon)|R|.$
\end{fact}
\begin{proof}
  Suppose the claim is false and let $V_i, i \in [\ell]$ be a cluster satisfying $\deg_R(V_i) < (c - 2d - 2\varepsilon)|R|.$ Then in the pure graph $G',$ we have:
  \begin{align*}
  \sum_{v \in V_i} \deg_{G'}(v) &\le \sum_{j \colon V_iV_j \in E(R)} e_{G'}(V_i, V_j) + e_{G'}(V_i, V_0) \le \deg_R(V_i) \frac{n}{|R|}|V_i| + |V_i||V_0| \\&< (c - 2d - 2\varepsilon) |V_i|n + \varepsilon |V_i| n = (c - 2d - \varepsilon) |V_i| n.
  \end{align*}
  This contradicts $\delta(G) \ge cn$ and property~\ref{reg_lemma_min_degree} given by Lemma~\ref{reg_lemma}.
\end{proof}

The so-called Slicing Lemma states that large subsets of a regular pair are also regular with slightly worse parameters.
\begin{lemma}[Slicing Lemma \cite{komlos1996szemeredi}]
	\label{slicing_lemma}
	Let $(A, B)$ be an $(\varepsilon, d)$-regular pair with density $d,$ and, for some $\alpha > \varepsilon,$ let $A' \subseteq A, \, |A|' \geq \alpha |A|, \, B' \subseteq B, \, |B'| \geq \alpha |B|.$ Then $(A', B')$ is an $\varepsilon'$-regular pair with $\varepsilon' = max\{\varepsilon / \alpha, 2 \varepsilon\}$ and for its density $d'$ we have $|d' - d| < \varepsilon.$
\end{lemma}

We also need the incredibly useful result of Koml{\'o}s, S{\'a}rk{\"o}zy and Szemer{\'e}di, known as the \textit{Blow-up} lemma, which states that $(\varepsilon, d)$-super-regular pairs behave like complete bipartite graphs in terms of containing subgraphs of bounded degree.
\begin{lemma}[Blow-up lemma \cite{komlos1997blow}]
	\label{blowup_lemma}
  Given a graph $R$ of order $r$ and positive parameters $\delta, \Delta,$ there exists a positive $\varepsilon = \varepsilon(\delta, \Delta, r)$ such that the following holds. Let $n_1, n_2, \dots, n_r$ be arbitrary positive integers and let us replace the vertices $v_1, v_2, \dots, v_r$ of $R$ with pairwise disjoint sets $V_1, V_2, \dots, V_r$ of sizes $n_1, n_2, \dots, n_r$ (blowing up). We construct two graphs on the same vertex set $V = \cup V_i.$ The first graph $F$ is obtained by replacing each edge $\{v_i, v_j\}$ of $R$ with the complete bipartite graph between the corresponding vertex-sets $V_i$ and $V_j.$ A sparser graph $G$ is constructed by replacing each edge $\{v_i, v_j\}$ arbitrarily with an $(\varepsilon, \delta)$-super-regular pair between $V_i$ and $V_j.$ If a graph $H$ with $\Delta(H) \leq \Delta$ is embeddable into $F$ then it is already embeddable into $G.$
\end{lemma}
The following remark also appears in \cite{komlos1997blow}.
\begin{remark}
	When using the Blow-up Lemma, we usually need the following strengthened version: Given $c > 0,$ there are positive numbers $\varepsilon = \varepsilon(\delta, \Delta, r, c)$ and $\alpha = \alpha(\delta, \Delta, r, c)$ such that the Blow-up Lemma in the equal size case (all $|V_i|$ are the same) remains true if for every $i$ there are certain vertices $x$ to be embedded into $V_i$ whose images are \emph{a priori} restricted to certain sets $C_x \subseteq V_i$ provided that
\begin{enumerate}[label=(\roman*), noitemsep, topsep=-5pt]
	\item each $C_x$ within a $V_i$ is of size at least $c |V_i|,$
	\item the number of such restrictions within a $V_i$ is not more than $\alpha |V_i|.$
\end{enumerate}

\end{remark}

% !TeX root = main.tex
% !TEX root = main.tex

\section{Lower bound examples}
We present simple lower bound constructions showing our results are best possible in a certain sense. For $\eta = \eta(n) = o(1),$ the condition $\delta(G) \geq (1 - \frac{1}{r+1} + \eta)n$ when $r \geq 3$ or $\delta(G) \geq (3/4 + \eta)n$ when $r=2$, is not enough to guarantee an $r^{th}$ power of a Hamilton cycle with absolute discrepancy linear in $n.$ Moreover, for $\eta = 0,$ there exists a graph in which every $r^{th}$ power of a Hamilton cycle has discrepancy $0.$

First consider $r \geq 3.$ Let $t \ge 2$ be even and $V_1, \dots, V_{r+1}$ disjoint clusters of size $t.$ Additionally, let $V_0$ be a cluster of size $0 \le m \le t.$ We construct a graph $G$ on the vertex set $V = \stackrel{.}{\bigcup}_{i=0}^{r+1} V_i.$ We put an edge between any two vertices in different clusters and we put all edges connecting two vertices in $V_0.$ Let $n = |V| = (r+1)t + m$ and note that $\delta(G) = rt + m = \left(1 - \frac{1}{r+1} + \frac{m}{(r+1)((r+1)t + m)}\right)n.$

Next we describe the coloring $f$ of the edges. We color the edges incident to vertices in $V_0$ arbitrarily. For each $V_i, \; i \geq 1$ we denote half of its vertices as \textit{positive} and the other half as \textit{negative}. For a vertex $v \in V_i$ and any vertex $u \in V_j$ where $1 \leq j < i$ we set $f(v, u) = 1$ if $v$ is positive and $f(v, u) = -1$ if $v$ is negative.

If $G$ contains no $r^{th}$ power of a Hamilton cycle, there is nothing to prove. Otherwise, let $H^r$ be an arbitrary $r^{th}$ power of a Hamilton cycle in $G$ viewed as a $2r$-regular subgraph of $G.$ Call a vertex $v \in V \setminus V_0$ a \textit{bad} vertex if at least one of its neighbours in $H^r$ is in the cluster $V_0,$ otherwise call it \textit{good}. If a vertex $v \in V_i$ is good then in $H^r$ it has precisely two neighbours in each of the clusters $V_j, \; 1 \leq j \leq r+1, j \neq i.$ Note that for $i \geq 1$ at most $2$ vertices in $V_i$ can be adjacent to a vertex $v \in V_0,$ so there are at most $2m$ bad vertices in $V_i.$ Now consider only positive good vertices and their edges towards vertices in clusters with a smaller index. Thus, the number of edges labelled $1$ in $H^r$ is at least
\[ \sum_{i=1}^{r+1} 2(i-1)(t/2 - 2m) = r(r+1)(t/2-2m). \]
Hence, we have
\[ f(H^r) \geq -nr + 2r(r+1)(t/2-2m) \geq -5r(r+1)m. \]
Completely analogously, $f(H^r) \leq 5r(r+1)m.$ Therefore, when $m = 0,$ or equivalently $\eta = 0,$ we have $f(H^r) = 0.$ If $m = o(n),$ or equivalently $\eta = o(n),$ we get $|f(H^r)| = o(n).$

For $r=2,$ the following construction was given in \cite{balogh2020discrepancies}, where the case $r=1$ was considered. Let $G$ be the $4$-partite Tur\'an graph on $n=4k$ vertices, so $\delta(G) = \frac{3}{4}n.$ Color all edges incident to one of the parts with $-1$ and the rest with $1.$ Any square of a Hamilton cycle contains $4k$ edges labelled $-1$, exactly $4$ for each vertex in the special class. As it has a total of $8k$ edges, its discrepancy is $0.$ Similarly as above, we can add $m = o(n)$ vertices connected to every other vertex and still any square of a Hamilton cycle has absolute discrepancy $o(n).$

% !TeX root = main.tex
% !TEX root = main.tex

\section{Outline of the proofs of Theorems \ref{main_thm} and \ref{r2_thm}}
Our proof follows a very similar structure to that of Balogh, Csaba, Pluh\'ar and Treglown \cite{balogh2020discrepancy} for $K_r$-tilings. 

We start by applying the regularity lemma on $G$ to obtain the reduced graph $R$ and the corresponding edge labelling $f_R.$ 

Before proving the P\'{o}sa-Seymour conjecture for large $n,$ Koml\'os, S\'ark\"ozy and Szemer\'edi \cite{approxposa}, proved an approximate version, namely they proved it for $n$-vertex graphs with minimum degree at least $(1 - 1/(r+1) + \varepsilon)n.$ We make slight modifications to their proof to establish two important claims.

We prove that a $K_{r+1}$-tiling of $R$ with linear discrepancy with respect to $f_R$ can be used to construct an $r^{th}$ power of a Hamilton cycle in $G$ with linear discrepancy with respect to $f.$ Combined with Theorem \ref{disc_hajnal_thm}, this is enough to deduce the case $r=2$ (Theorem \ref{r2_thm}).

Next, we prove a very useful property of $R$: suppose that $F$ is a small subgraph of $R$ and there are two $C^r$-templates of $F$ covering the vertices the same number of times, but having different discrepancies with respect to $f_R.$ Then in $G$ there exists the $r^{th}$ power of a Hamilton cycle with linear discrepancy with respect to $f.$

We assume $G$ has no $r^{th}$ power of a Hamilton cycle with large absolute discrepancy. From this point on, we only `work' on the reduced graph $R.$ To use the last claim, we need a subgraph $F$ on which we can find two different $C^r$-templates, so the simplest subgraph we can study is an $(r+2)$-clique. We prove that every $(r+2)$-clique in $R$ is either a copy of $K_{r+2}^+, K_{r+2}^-, (K_{r+2}, +)$-star or $(K_{r+2}, -)$-star. As $R$ has large minimum degree, every clique of size $k \leq r+2$ can be extended to a clique of size $r+2.$ This shows that every clique of size $k \leq r+2$ is either a copy of $K_k^+, K_k^-, (K_k, +)$-star or $(K_k, -)$-star.

By the Hajnal-Szemer\'edi theorem, we can find a $K_{r+1}$-tiling $\calT$ of $R.$ The previous arguments show that only four types of cliques appear in $\calT$ and $\calT$ has a small discrepancy with respect to $f_R.$ This tells us that the numbers of each of the four types of cliques in $\calT$ are balanced in some way. 

We consider two cliques in $\calT$ of different types. For several relevant cases when there are many edges between the two cliques, we construct two $C^r$-templates of different discrepancies, which contradicts our assumption by the claim about $C^r$-templates. Finally, this restricts the number of edges between different cliques which leads to a contradiction with the minimum degree assumption on $R.$

% !TeX root = main.tex
% !TEX root = main.tex

\section{Using $C^r$-tilings of $R$}

In \cite{approxposa} Koml\'os, S\'ark\"ozy and Szemer\'edi proved an approximate
version of the P\'osa-Seymour conjecture. More precisely, they show that for any
$\eta > 0,$ a graph $G$ on $n$ vertices, for sufficiently large $n$, with minimum
degree at least $(r/(r+1)+\eta)n$ contains the $r^{th}$ power of a Hamilton
cycle. Their argument starts by applying the regularity lemma to get an
$(\varepsilon,\eta/3)$-regular partition $(V_i)_{i = 0}^{t}$ of $V(G)$ and the
Hajnal-Szemer\'edi theorem to obtain a $K_{r+1}$-tiling of its
$(\varepsilon,\eta/3)$-reduced graph $R$. Let $K_1, \dotsc, K_s$ be the $(r+1)$-cliques
in this tiling. Then they proceed to find short paths $\calP_1, \dotsc, \calP_s$ whose $r^{th}$ powers are in $G$,
each $\calP_i$ `connecting' subsequent cliques $K_i, K_{i+1}$ (where we denote $K_{s+1} = K_1$) in the tiling, and
`attach' the exceptional vertices---the ones in $V_0$ together with some other
vertices not respecting a certain degree condition---to these paths, thus obtaining an $(\varepsilon, \eta/9)$-super-regular partition. Finally,
for each clique $K_i$, the Blow-up Lemma~\cite{komlos1997blow} is applied to find the
$r^{th}$ power of a Hamilton path on the set of unused vertices in $G$
corresponding to this clique, which together with paths $\calP_1, \dotsc, \calP_s$ closes
the $r^{th}$ power of a Hamilton cycle.

We need a slightly more general result for our application. Instead of a
$K_{r+1}$-tiling in $R$, we assume a $C^r$-tiling $\calT$ of $R$ is given in which most of
the cycles $C_i \in \calT$ have length $r+1$. Then we proceed similarly as
outlined above, where for each $C_i \in \calT$ we choose an $(r+1)$-clique in $C_i^r$ 
to represent $K_i$ for which the connecting path $\calP_i$ is constructed. In order
to argue about the discrepancy of the found $r^{th}$ power of a Hamilton cycle
$H^r$ in $G$, we explicitly state a property of the construction in
\cite{approxposa}, which is that most of the edges in $H^r$ come from the given
tiling and that these edges are used in a balanced way. Intuitively, if the
$C^r$-tiling of $R$ is given together with a function $f_R$ as in \eqref{deffR},
then the discrepancy of the found $H^r$ in $G$ is at most $\alpha n$ away from
$m \cdot f_R(\calT),$ for an arbitrarily
small $\alpha > 0$, where we remind the reader that $m$ is the size of each cluster. Additionally, we show that, given two similar $C^r$-tilings,
we can find similar $r^{th}$ powers of a Hamilton cycle. We make these notions
precise in the following statement. 

\begin{proposition}
  \label{posa_seymour_prop}
  For any integer $r \geq 2$ and any $\alpha, \eta > 0$ there exist $\ell_0, n_0 \in
  \naturals$ and $\varepsilon > 0$ such that the
  following holds. Suppose $G$ is a graph on $n \ge n_0$ vertices with $\delta(G) \geq
  (r/(r+1)+\eta)n$ and let $f \colon E(G) \rightarrow \{-1, 1\}$ be its edge
  labelling. Let $(V_i)_{i=0}^{\ell},$ where $\ell \geq \ell_0,$ be an
  $(\varepsilon,\eta/3)$-regular partition of $G$ with respect to $f$ with $|V_i| = m, \forall i \in [\ell],$ and let $R$ its
  $(\varepsilon,\eta/3)$-reduced graph with $f_R \colon E(R) \to \{-1,1\}$ as defined
  in \eqref{deffR}. Suppose we are given two $C^r$-tilings $\mathcal{T}_1 =
  \mathcal{K} \stackrel{.}{\cup} \mathcal{C}_1$ and $\mathcal{T}_2 = \mathcal{K} \stackrel{.}{\cup}
  \mathcal{C}_2$ of $R$ such that
  \begin{itemize}
    \item $\mathcal{K}$ consists only of $(r+1)$-cycles, and
    \item $|\mathcal{C}_i| \leq 10r^2$, for $i \in \{1,2\}$.
  \end{itemize}
  Then there exist $r^{th}$ powers of Hamilton cycles $H_1^r, H_2^r \subseteq G$
  such that
  \begin{enumerate}[label=(\roman*)]
    \item\label{discofHr} $|f(H^r_i) - m f_R(\mathcal{T}_i)| \leq \alpha n$, for
      $i \in \{1, 2\}$, and
    \item\label{disc_diff} $|f(H_1^r) - f(H_2^r)| \geq \frac{m}{2}
      |f_R(\mathcal{C}_1) - f_R(\mathcal{C}_2)|.$
  \end{enumerate}
\end{proposition}

\begin{proof}
  %\MT{Maybe:
  The proof is almost a one-to-one copy of the argument by Koml\'os, 
    S\'ark\"ozy and Szemer\'edi in \cite{approxposa}.  Nevertheless, for completeness we present the full proof which follows the rough outline given in the paragraphs above the statement. Several claims with lengthy proofs which do not require any adaptation from \cite{approxposa} to our setting are presented in the appendix.

  We may assume $1/n \ll 1/\ell_0 \ll \varepsilon \ll \alpha \ll \eta \ll 1/r.$
  Recall that $|V_0| \leq \varepsilon n$ and let $m = |V_i|$ for $1 \leq i \leq \ell.$  Let $\mathcal{K} = \{C_1^r, \dots, C_s^r\}$ and $\mathcal{C}_1 = \{C_{s+1}^r, \dots, C_{s+q}^r\}$ where $s \leq \frac{\ell}{r+1}$ and $q \leq 10r^2.$ The cycles $C_1, \dots, C_s$ are of size $r+1,$ that is, $C_1^r, \dots, C_s^r$ are $(r+1)$-cliques where each edge has multiplicity $2.$
  For $i \in [s+q],$ let $V^i_1, V^i_2, \dots, V^i_{r+1}$ be the first $r+1$ vertices of $C_i$ and let $K_i$ denote the clique on these $r+1$ vertices. Let us denote $d = \eta / 3.$ From Fact~\ref{degRfact}, we have that for $\eta' = \eta/4,$
\begin{equation} \label{eq:degR}
  \delta(R) \ge (1 - 1/(r+1) + \eta') \ell,
\end{equation}
and consequently,
\begin{equation} \label{eq:common_neigh}
  \forall V_1, \dots, V_{r+1} \in V(R), \; \deg_R(\{V_1, \dots, V_{r+1}\}) \ge (r+1)\eta' \ell.
\end{equation}
  The first step of the proof is to find short connecting paths $\mathcal{P}_i$ between cliques $K_i$ and $K_{i+1}$ for all $i \in [s+q],$ where we denote $K_{s+q+1} = K_1.$ 

\begin{restatable}{claim}{claimconnectingpaths} \label{cl:conn_paths}
In $G$ there exist vertex disjoint $r^{th}$ powers of paths $\calP_1^r, \dots, \calP^r_{s+q},$ where $\calP_i = (p^i_1, p^i_2, \dots, p^i_{t_i})$ such that the following holds for all $i \in [s+q]$:
\begin{enumerate}[label=P\arabic*)]
  \item \label{enum:connects} $\calP_i^r$ connects the cliques $K_i$ and $K_{i+1},$ that is, for any $j \in [r+1],$ $p^i_j \in V^i_j$ and $p^i_{t_i-r-1+j} \in V^{i+1}_j$ (where we identify $V^{s+q+1}_j = V^1_j$);
  \item \label{enum:extensible_endpoints} for any $j \in [r],$
  \[ |N(\{p^i_1, p^i_2, \dots, p^i_j\}) \cap V^i_{j+1}|, \, |N(\{p^i_{t_i-j+1}, p^i_{t_i-j+2}, \dots, p^i_{t_i}\}) \cap V^{i+1}_{r+1-j}| > (d-\varepsilon)^{2r} m/4; \]
  \item $t_i \leq O(r^3).$
\end{enumerate}
Moreover, these paths can be constructed one by one such that the paths $\calP_1, \dots, \calP_{s-1}$ only depend on $\calK.$
\end{restatable}

Condition \ref{enum:extensible_endpoints} will ensure we can later close the $r^{th}$ power of a Hamilton cycle inside each clique. In order to prove part \ref{disc_diff} of the statement, we will use the fact we can construct the connecting paths such that the first $s-1$ of them only depend on $\calK.$ Roughly speaking, this will allow us to construct two $r^{th}$ powers of Hamilton cycles which only differ in a few edges apart from the edges which touch the vertices in $\calC_1$ (or equivalently in $\calC_2).$ This way, the difference of the discrepancies of the two $r^{th}$ powers of Hamilton cycles is essentially controlled by $f_R(\calC_1) - f_R(\calC_2),$ as desired. The proof of Claim~\ref{cl:conn_paths} is nearly a direct copy of an argument from \cite{approxposa} so we defer it to the appendix.

  Next we add some more vertices to the exceptional set $V_0$ with the goal of obtaining an $(\varepsilon, d/3)$-super-regular partition, which will be verified later.
  From a cluster $V_j^i$ in a cycle $C_i$ we move to $V_0$ all vertices $v$ not used on the paths $\mathcal{P}_1, \dots, \mathcal{P}_{s+q}$ for which there is a $j'$ such that
\[ \{V_j^i, V_{j'}^i\} \in E(C_i^r) \text{ and } \deg_{f_R(V_j^i, V_{j'}^i)}(v, V_{j'}^i) < (d-\varepsilon) |V_{j'}^i|, \]
where we consider the original clusters $V_{j'}^i,$ that is, we also consider the vertices already used in the connecting paths in this calculation. Because of $(\varepsilon, d)$-regularity and because $|C_i| \leq 10r^2$ for all $i \in [s+q]$ by definition of a $C^r$-tiling, there are at most $|E(C_i^r)|\cdot \varepsilon m \leq 10 r^3 \varepsilon m$ such vertices in each cluster $V_j^i.$ Then we move the smallest possible number of vertices to $V_0$ so that each cluster has the same number of vertices. We still write $V_0$ for the enlarged exceptional set which now satisfies $|V_0| \leq 11r^3 \varepsilon n$ and we use $m'$ to denote the new number of vertices in each cluster. Recall that we have used at most $(s+q) \cdot O(r^3) = O(\ell r^2)$ vertices for the connecting paths. Hence, 
\begin{equation}\label{eq:m_prime}
  m' = \frac{n - |V_0| - O(\ell r^2)}{\ell} \ge \frac{n}{\ell} \cdot (1 - 12r^3 \varepsilon) \ge (1 - 12r^3 \varepsilon)m.
\end{equation}
For each vertex $v \in V_0$ we find all indices $i \in [s-1]$ such that
  \begin{equation} \label{eq:canbeassigned}
    \deg(v, V_j^{i+1}) \geq (d-\varepsilon) |V_j^{i+1}|, \; \forall j \in [r+1].
  \end{equation}
    Again, in the above we consider the original clusters $V_j^{i+1}.$ Let $x$ denote the number of such indices $i$. Note that $v$ has at most $\varepsilon n$ edges to the initial exceptional set and for every $i$ for which \eqref{eq:canbeassigned} is not satisfied at most $(r + d - \varepsilon)m$ edges to the clusters in the clique $K_{i+1}.$ Therefore, we have
\[ \left(1 - \frac{1}{r+1} + \eta\right)n \leq \deg(v) \leq \varepsilon n + \left(|C_1| + \left(\sum_{i=s+1}^{s+q} |C_i|\right)\right)m + x(r+1)m + (s - 1 - x)(r + d - \varepsilon) m. \]
  Using $q \leq 10r^2, |C_i| \leq 10r^2,$ that $\varepsilon$ is sufficiently small and $\ell$ is sufficiently large, a simple calculation yields $x \geq \eta \ell / 2.$ We assign $v$ to one of these cliques such that no clique is assigned too many vertices. That is, we partition $V_0$ into small sets $A_1, \dots, A_{s-1}$ such that every $v \in A_i$ satisfies \eqref{eq:canbeassigned} for the index $i.$ We partition $V_0$ by considering vertices one by one. The next vertex $v \in V_0$ is put into the currently smallest set $A_i$ among the indices $i$ for which \eqref{eq:canbeassigned} holds with respect to $v.$ It is easy to see that in the end the sets $|A_i|$ have size at most $\frac{|V_0|}{\eta \ell /2} \leq \varepsilon_1 m',$ where $\varepsilon_1 = \frac{100 r^3 \varepsilon}{\eta}.$
\begin{restatable}{claim}{claimaddingvertices} \label{cl:adding_vertices_to_paths}
The $r^{th}$ power of a path $\mathcal{P}_i$ can be extended to contain all vertices in $A_i$ by using additional three vertices from each of the clusters $V^{i+1}_j, \; j \in [r+1]$ per added vertex. This can be done so that property \ref{enum:extensible_endpoints} holds with respect to the first and last $r$ vertices of the new path. 
\end{restatable}
Again, we use the argument from \cite{approxposa} and defer the proof of this claim to the appendix.

For any $i \in [s+q], j \in [|C_i|],$ let $W_j^i \subseteq V_j^i$ be the set of unused vertices in $V_j^i.$ By construction, for each $i \in [s+q],$ the sets $W_j^i$ are of the same size for all $j \in [|C_i|].$ From $V_j^i,$ we have used at most $\ell \cdot O(r^3)$ vertices to form the initial connecting paths and then exactly $3 |A_i| \le 3\varepsilon_1 m'$ vertices to extend the path to include the vertices in $A_i.$ Hence,
  \begin{equation}
    \label{enough_left}
    |W_j^i| \ge m' - \ell \cdot O(r^3) - 3 \varepsilon_1m' \ge (1 - \frac{300r^3 \varepsilon}{\eta})(1 - 12r^3\varepsilon)m - \ell \cdot O(r^3) \ge (1 - \frac{\alpha}{6r})m \eqqcolon m'',
  \end{equation}
    
  where we chose $\varepsilon$ to be small enough and $n_0$ to be large enough.

  Let us now verify that whenever $\{V_j^i, V_{j'}^i\} \in E(C_i^r),$ the bipartite subgraph $G_\sigma[W_j^i, W_{j'}^i]$ is $(\varepsilon, d/3)$-super-regular, where $\sigma = f_R(V_j^i, V_{j'}^i).$ Let $i, j, j'$ be such that $\{V_j^i, V_{j'}^i\} \in E(C_i^r).$ Recall that for any $v \in W_j^i,$ we have $\deg_\sigma(v, V_{j'}^i) \ge (d-\varepsilon)m$ as otherwise $v$ would have been moved to $V_0$ and so it would not be in $W_j^i.$ Thus,
  \[ \deg_\sigma(v, W_{j'}^i) \ge (d-\varepsilon)m - (|V_{j'}^i| - |W_{j'}^i|) \ge (d-\varepsilon)m - (m-m'') \ge (d - \varepsilon - \frac{\alpha}{6r})m \ge \frac{d}{3}m \ge \frac{d}{3} |W_{j'}^i|, \]
  where we chose $\alpha$ small enough compared to $\eta.$ Now, consider subsets $X \subseteq W_j^i, \, Y \subseteq W_{j'}^i$ such that $|X| > \varepsilon |W_j^i|$ and $|Y| > \varepsilon |W_{j'}^i|.$ Because $G_{\sigma}[V_j^i, V_{j'}^i]$ is $(\varepsilon, d)$-regular, it follows that $d_{G_\sigma}(X, Y) \ge d-\varepsilon > d/3,$ so the bipartite graph $G_\sigma[W_j^i, W_{j'}^i]$ is indeed $(\varepsilon, d/3)$-super-regular.

  Finally, the endpoints of the connecting paths satisfy \ref{enum:extensible_endpoints} so we can apply Lemma \ref{blowup_lemma} and the subsequent remark to close the $r^{th}$ power of a Hamilton cycle $H^r$ such that between clusters $V_x, V_y$ we only use edges of color $f_R(V_x, V_y).$

  From \eqref{enough_left}, it easily follows that $H^r$ satisfies \ref{discofHr}. Indeed, for every edge $(V_x, V_y)$ in some $C_i^r, \; i \in [s+q],$ we included in $H^r$ at least $\mul_{C_i^r}(V_x, V_y) \cdot (1 - \frac{\alpha}{6r})m$ corresponding edges of color $f_R(V_x, V_y).$ Apart from these, there are at most $(|V_0| + \frac{\alpha}{6r}n) \cdot 2r \leq \frac{\alpha}{2} n$ edges in $H^r.$ Hence,
  \[ |f(H^r) - m f_R(\calT_1)| \leq \frac{\alpha}{6r}m |f_R(\calT_1)| + \frac{\alpha}{2} n \leq \frac{\alpha}{6r} m \cdot 2r \ell + \frac{\alpha}{2} n < \alpha n. \]

  Now suppose we are given two tilings $\calT_1, \calT_2$ as in the statement of the proposition. We construct two $r^{th}$ powers of a Hamilton cycle $H_1^r$ and $H_2^r$ with a small modification to the above procedure to ensure they do not differ too much. We think of the above algorithm as three stages of adding edges to a subgraph which eventually becomes the $r^{th}$ power of a Hamilton cycle. We use $\calI_1$ and $\calI_2$ to denote the runs of the algorithm on $\calT = \calT_1$ and $\calT = \calT_2,$ respectively.

  The first stage of connecting the cliques is done exactly as above. Note that the paths $\mathcal{P}_1, \dots, \mathcal{P}_{s-1}$ are the same in both instances. Let $V_0^1$ denote the exceptional set and $U^1$ the set of vertices used for the connecting paths in $\calI_1.$ Analogously define $V_0^2$ and $U^2$ with respect to $\calI_2.$ Note that $|U^1 \Delta U^2| \leq 2(q + 1) \cdot O(r^3) = O(r^5)$ because the two runs only deviate after finding the first $s-1$ paths and each path has length $O(r^3).$ Since the degree of any vertex in the $r^{th}$ power of a Hamilton cycle is $2r,$ this implies the two subgraphs differ in $O(r^6)$ edges at this point.

  We can assume that the new number of vertices per cluster $m'$ is the same in $\calI_1$ and $\calI_2,$ otherwise simply add a few vertices to $V_0^j$ if $m'$ is larger in $\calI_j.$ Set $V_0 = V_0^1 \cup V_0^2$ and note that now $|V_0| \leq 22r^3 \varepsilon n$ which does not affect the above argument. To make the two found subgraphs similar, we treat $V_0 \setminus U^1$ and $V_0 \setminus U^2$ as the new exceptional sets for $\calI_1$ and $\calI_2,$ respectively. In the second stage, we add the exceptional vertices to paths $\mathcal{P}_1, \dots, \mathcal{P}_{s-1}.$ As noted above, the cliques to which a vertex can be assigned do not depend on the connecting paths. Therefore, we can assign the vertices in $V_0 \setminus (U^1 \cup U^2)$ to the same cliques for the two runs. We can then embed the vertices from $V_0 \setminus (U^1 \cup U^2)$ in exactly the same way in $\calI_1$ and $\calI_2.$ Thus $\calI_1$ and $\calI_2$ only deviate after embedding all but at most $|U^1 \Delta U^2| = O(r^5)$ vertices. Each added vertex extends the path by $O(r)$ vertices, so this stage introduces at most $O(r^7)$ edges to $H_1^r \Delta H_2^r.$

  Finally, in the third stage, we find connecting paths inside cliques $K_i.$ The edges in these paths have the same labels in $\calI_1$ and $\calI_2$ for all cliques $K_i, \; i \in [s].$ The number of vertices left in cliques $K_1, \dots, K_s$ might differ in $\calI_1$ and $\calI_2$ at this point, but in total only by $|U^1 \Delta U^2| \cdot O(r).$ Hence, the difference between the values of edges in cliques $K_1, \dots, K_s,$ introduced at this stage is at most $O(r^7).$ If $f_R(\calC_1) = f_R(\calC_2),$ there is nothing to prove so let us assume $f_R(\calC_1) \neq f_R(\calC_2).$ Since for each edge in $\mathcal{C}_1$ we add $m' \geq (1 - \frac{\alpha}{6r})m$ edges to $H_1^r$ and analogously for $\mathcal{C}_2$ and $H_2^r,$ we have
  \[ |f(H_1^r) - f(H_2^r)| \geq \left(1 - \frac{\alpha}{6r}\right)m \cdot |f_R(\mathcal{C}_1) - f_R(\mathcal{C}_2)| -  O(r^7) \geq \frac{m}{2} \big|f_R(\mathcal{C}_1) - f_R(\mathcal{C}_2)\big|, \]
  where we used $|f_R(\calC_1) - f_R(\calC_2)| \ge 1$ and $m \gg r.$
\end{proof}

\section{Proofs of Theorems \ref{main_thm} and \ref{r2_thm}}
% !TeX root = main.tex
% !TEX root = main.tex

It is enough to prove the theorems for $\eta \ll 1/r.$ We define additional constants $\alpha, \beta, \gamma, \varepsilon, d > 0$ and $n_0, \ell_0, L_0 \in \naturals$ such that 
\[ 0 < 1/n_0 \ll \gamma \ll 1/L_0 \leq 1/\ell_0 \ll \varepsilon \ll \alpha \ll \beta \ll d = \eta/3 \ll 1/r, \]
where $\ell_0, n_0, \varepsilon$ are chosen so that Proposition~\ref{posa_seymour_prop} holds for $r, \alpha, \eta$ and $L_0$ is the constant obtained from Lemma \ref{reg_lemma} with parameters $\varepsilon/2, \ell_0.$ Let $G$ be a graph with $n \geq n_0$ vertices and $f \colon E(G) \rightarrow \{-1, 1\}$ an edge labelling as in the theorem statements.

Both proofs start by applying Lemma \ref{reg_lemma} to $G$ with parameters $\varepsilon/2, d$ and $\ell_0$ and the subsequent remark to make the number of clusters divisible by $r+1.$ We thus obtain an $(\varepsilon, d)$-regular partition of $G$ with respect to $f$ and the reduced graph $R$ whose vertices are clusters $V_1, V_2, \dots, V_\ell,$ where $r+1 \mid \ell$ and each of the clusters $V_i$ is of size $m.$ We also have the exceptional set $V_0$ of size at most $\varepsilon n.$ The reduced graph $R$ inherits the edge labelling $f_R$ as given in \eqref{deffR}. 

The following claim is a direct consequence of Proposition~\ref{posa_seymour_prop}~part~\ref{discofHr}.
\begin{claim}
  \label{tilingtoHr_claim}
  Suppose there exists a $K_{r+1}$-tiling $\calT$ of $R$ such that $|f_R(\calT)| \geq \beta \ell.$ Then, in $G$ there exists the $r^{th}$ power of a Hamilton cycle $H^r$ satisfying $|f(H^r)| \geq \gamma n.$ \qed
\end{claim}
\begin{proof}
  Suppose we are given a $K_{r+1}$-tiling $\calT$ of $R$ satisfying $|f_R(\calT)| \ge \beta \ell.$ Recall that we can view $\calT$ as a $C^{r}$-tiling of $R$ so by Proposition~\ref{posa_seymour_prop}~part~\ref{discofHr}, in $G$ there exists an $r^{th}$ power of a Hamiton cycle $H^r$ satisfying
  \[ |f(H^r) - m f_R(\calT)| \le \alpha n, \]
  which implies
  \[ |f(H^r)| \ge m \beta \ell - \alpha n \ge (1 - \varepsilon) \beta n - \alpha n > \gamma n, \]  
  where we chose $\alpha, \gamma, \varepsilon$ to be small enough compared to $\beta.$
\end{proof}

Now we resolve the case $r=2$ which can be easily deduced from Theorem \ref{disc_hajnal_thm} and Claim \ref{tilingtoHr_claim}.
\begin{proof}[Proof of Theorem \ref{r2_thm}]
  Recall that $\delta(G) \geq (3/4 + \eta)n,$ so by Fact~\ref{degRfact}, we get $\delta(R) \geq (3/4 + \eta/4) |R|.$
	Let $\beta$ be the value of $\gamma$ given by Theorem \ref{disc_hajnal_thm} with parameters $r = 3$ and $\eta / 4.$	Applying Theorem \ref{disc_hajnal_thm} to the reduced graph $R$ we obtain a $K_3$-tiling $\calT$ of $R$ of absolute discrepancy at least $\beta \ell.$ By Claim \ref{tilingtoHr_claim}, $G$ contains the square of a Hamilton cycle with absolute discrepancy at least $\gamma n.$
\end{proof}

In the rest of the paper we prove Theorem~\ref{main_thm}. Recall that $\delta(G) \ge (1 - 1/(r+1) + \eta)n,$ so by Fact~\ref{degRfact}, the reduced graph satisfies
\begin{equation}
  \label{degreeofR}
  \delta(R) \geq \left(1 - \frac{1}{r+1} + \frac{\eta}{4}\right) \ell,
\end{equation}
where we remind the reader that $\ell = |R|.$ The following simple observation follows directly.
\begin{claim}
  \label{common_neigh}
  For any $v_1, v_2, \dots, v_{r+1} \in R,$ we have $|N(\{v_1, v_2, \dots, v_{r+1}\})| \geq (r+1) \eta \ell/4.$ \qed
\end{claim}

Next we derive an additional claim from Proposition~\ref{posa_seymour_prop}.
\begin{claim}
  \label{template_claim}
  Let $F$ be a subgraph of $R$ on at most $10r$ vertices and let $\calF_1 = \{C_{11}, C_{12}, \dots, C_{1s_1}\}$ and $\calF_2 = \{C_{21}, C_{22}, \dots, C_{2s_2}\}$ be two $C^r$-templates of $F$ such that each vertex of $F$ appears exactly $k$ times in $\calF_1$ and $k$ times in $\calF_2$ for some $k \leq 10r.$ If $\calF_1$ and $\calF_2$ have different discrepancies with respect to $f_R$, then in $G$ there exists the $r^{th}$ power of a Hamilton cycle $H^r$ satisfying $|f(H^r)| \geq \gamma n.$
\end{claim}
\begin{proof}
  Split each of the clusters $V_i, 1 \le i \le \ell,$ of the regular partition into $k$ clusters $V_{i, 1}, \dots, V_{i, k}$ of size $m' = \lfloor m/k\rfloor$ and put the remaining vertices in $V_0.$ Let $c(V_i) = \{V_{i, 1}, \dots, V_{i, k}\}.$ Define $R'$ as a blow-up of $R$ on the vertex set $\bigcup_{i=1}^\ell c(V_i)$ with the edge coloring $f_{R'}$ in the natural way: if $V_i$ and $V_j$ were joined by an edge in $R,$ then put a complete bipartite graph between $c(V_i)$ and $c(V_j)$ with all edges of color $f_R(V_i, V_j).$ Let $F'$ denote the corresponding blow-up of $F.$ Formally, we define $V(F') = \bigcup_{V_i \in V(F)} c(V_i)$ and $E(F') = \bigcup_{(V_i, V_j) \in E(F)} \{ XY \; \vert \; X \in c(V_i), Y \in c(V_j) \}.$ By Lemma \ref{slicing_lemma}, we get that every pair of clusters joined by an edge in $R'$ form an $\varepsilon'$ regular pair in $G$ for $\varepsilon' = 2k \varepsilon$ and have density at least $d' = d - \varepsilon.$ 
  
  By Remark~\ref{rem_divisibility}, we may assume that the number of vertices of $R' \setminus F'$ is divisible by $r+1$ by moving at most $r$ clusters into the exceptional set $V_0'$ which now satisfies $|V_0'| \le \varepsilon' n.$ Note that $|R'| = k|R|$ and $\delta(R') = k \delta(R),$ implying $\delta(R') \geq (1 - \frac{1}{r+1} + \frac{\eta}{4}) |R'|$ and therefore, $\delta(R' \setminus F') \geq (1 - \frac{1}{r+1})|R'|.$ Applying Theorem \ref{hajnal_szemeredi}, we obtain a $K_{r+1}$-tiling $\calK$ of $R' \setminus F'.$ Using $\calF_1$ we construct a $C^r$-tiling $\calC_1$ of $F'$ as follows. In the cycles in $\calF_1,$ simply replace every occurrence of a vertex $V_i$ with a different vertex from $c(V_i).$ It is easy to verify that $\calC_1$ is a $C^r$-tiling of $F'.$ Analogously, we construct $\calC_2$ from $\calF_2.$ 
  
  Let $\calT_1 = \calK \stackrel{.}{\cup} \calC_1$ and $\calT_2 = \calK \stackrel{.}{\cup} \calC_2.$ We apply Proposition \ref{posa_seymour_prop} to obtain two $r^{th}$ powers of Hamilton cycles $H_1^r$ and $H_2^r$ which satisfy \ref{disc_diff} with respect to $m'$ and $f_{R'}$. Finally, note that $f_{R'}(\calT_i) = f_{R'}(\calK) + f_{R'}(\calC_i)$ and $f_{R'}(\calC_i) = f_R(\calF_i)$ for $i=1,2.$ Hence, we have
  \[ |f(H_1^r) - f(H_2^r)| \geq \frac{m'}{2} |f_{R'}(\calT_1) - f_{R'}(\calT_2)| \geq \frac{m}{2(k+1)} \cdot 1 \geq 2 \gamma n. \]
  Therefore, at least one of $H_1^r, H_2^r$ has absolute discrepancy at least $\gamma n.$
\end{proof}

% !TeX root = main.tex
% !TEX root = main.tex

First we resolve the case $r=2$ which can be easily deduced from Theorems \ref{disc_hajnal_thm} and Claim \ref{tilingtoHr_claim}.
\begin{proof}[Proof of Theorem \ref{r2_thm}]
  Recall that $\delta(G) \geq (3/4 + \eta)n,$ so by Fact \ref{degRfact}, we get $\delta(R) \geq (3/4 + \eta/4) |R|.$
	Let $\beta$ be the value of $\gamma$ given by Theorem \ref{disc_hajnal_thm} with parameters $r = 3$ and $\eta / 4.$	Applying Theorem \ref{disc_hajnal_thm} to the reduced graph $R$ we obtain a $K_3$-tiling $\calT$ of $R$ of absolute discrepancy at least $\beta \ell.$ By Claim \ref{tilingtoHr_claim}, $G$ contains the square of a Hamilton cycle with absolute discrepancy at least $\gamma n.$
\end{proof}

% !TeX root = main.tex
% !TEX root = main.tex

In the rest of the paper, we finish the proof of Theorem \ref{main_thm}.
\begin{proof}[Proof of Theorem \ref{main_thm}]
  We prove the theorem by contradiction, so we assume that $G$ does not contain the $r^{th}$ power of a Hamilton cycle with absolute discrepancy at least $\gamma n.$  Recall that we have $r \geq 3,\, \delta(G) \geq \left(1 - \frac{1}{r+1} + \eta\right) n$ and $\delta(R) \ge (1 - \frac{1}{r+1} + \eta/4)\ell.$

  The following claim shows that $R$ is highly structured with respect to $f_R.$
  \begin{claim}
    \label{cliques_cor}
    Let $K$ be a clique in $R$ of size $k \leq r+2.$ Then $K$ is a copy of one of the following: $K_k^+, K_k^-, (K_k, +)$-star or $(K_k, -)$-star.
  \end{claim}
  \begin{proof}
    First we prove the claim for $(r+2)$-cliques. Let $K = \{v_1, v_2, \dots, v_{r+2}\}$ be an $(r+2)$-clique in $R.$ Define $C_1 = (v_1, v_2, v_3, v_4, v_5, \dots, v_{r+2})$ and $C_2 = (v_1, v_3, v_2, v_4, v_5, \dots, v_{r+2}).$ We can view $\{C_1\}$ and $\{C_2\}$ as $C^r$-templates on $K.$ Thus, by Claim \ref{template_claim}, we get $f_R(C_1^r) = f_R(C_2^r).$ Note that 
    \[ f_R(C_1^r) = 2\sum_{1 \leq i <j \leq r+2} f_R(v_i, v_j) - \sum_{i=1}^{r+2} f_R(v_i, v_{i+1}), \]
    where we denote $v_{r+3} = v_1.$ Hence, $0 = f_R(C_1^r) - f_R(C_2^r) = f_R(v_1, v_3) + f_R(v_2, v_4) - f_R(v_1, v_2) - f_R(v_3, v_4).$ As the enumeration of the vertices was arbitrary, for any distinct $a, b, c, d \in K,$ the following holds:
    \begin{equation}
      \label{square_equation}
      f_R(a, b) + f_R(c, d) = f_R(a, c) + f_R(b, d).
    \end{equation}
    
    The rest of the proof appears in \cite{balogh2020discrepancy}. We present a slightly shorter argument. Assume that $K$ is not monochromatic, so there exists a vertex $v \in K$ with $N_K^+(v), N_K^-(v) \neq \emptyset.$ Without loss of generality, assume $|N_K^+(v)| \geq 2$ and let $u \in N_K^-(v).$ Consider arbitrary distinct $x, y \in N_K^+(v).$ By \eqref{square_equation}, we get $f_R(x, v) + f_R(u, y) = f_R(x, y) + f_R(u, v).$ By, definition $f_R(x, v) = 1$ and $f_R(u, v) = -1,$ so this implies $f_R(u, y) = -1$ and $f_R(x, y) = 1.$ If $|N_K^-(v)| \geq 2,$ a completely analogous argument shows $f_R(u, y) = 1,$ a contradiction. From this we conclude $N_R^-(v) = \{u\}.$ Applying the same reasoning to every pair $x, y \in N_R^+(v),$ we get $f(x, y) = 1$ and $f(u, y) = -1$ for any $x,y \in N_R^+(v).$ In other words, $K$ is a $(K_{r+2}, -)$-star with $u$ as its head.

    Now, suppose $K$ is a clique in $R$ of size $k \leq r+2.$  By Claim \ref{common_neigh}, $K$ can be extended, vertex by vertex, to some clique $K'$ of size $r+2.$ The statement now easily follows from the result for $(r+2)$-cliques.
    \end{proof}

  By \eqref{degreeofR}, $R$ has large minimum degree, so we can apply Theorem \ref{hajnal_szemeredi} to obtain a $K_{r+1}$-tiling $\calT$ of $R.$ From Claim \ref{cliques_cor} we conclude there are only four types of cliques in $\calT.$ Let $A$ denote the set of $K_{r+1}^+$ in $\calT; \; B$ the set of $K_{r+1}^-$ in $\calT; \; C$ the set of $(K_{r+1}, +)$-stars in $\calT;$ and $D$ the set of $(K_{r+1}, -)$-stars in $\calT.$ 
	Under the assumption that $G$ does not have $r^{th}$ powers of Hamilton cycles with large discrepancy, we establish several claims about edges between cliques of different types. We state these claims here and defer their proofs to the end of the paper.

  \begin{restatable}{claim}{claimAtoB}
    \label{AtoB}
    Consider a vertex $x_1$ of a clique $X \in A$ and let $Y$ be a copy of $K^-_{r+1}$ in $B.$ Then we may assume $\deg(x_1, Y) \leq r-1.$
  \end{restatable}
 
  \begin{restatable}{claim}{claimAtoC}
    \label{AtoC}
    Consider a vertex $x_1$ of a clique $X \in A$ and let $Y$ be a $(K_{r+1}, +)$-star in $C.$ Then we may assume $\deg(x_1, Y) \leq r-1.$
  \end{restatable}

  \begin{restatable}{claim}{claimDtoC}
    \label{DtoC}
    Suppose $x_1$ is the head of a clique $X \in D$ and let $Y \in C.$ Then, we may assume $\deg(x_1, Y) \leq r-1.$
  \end{restatable}

  \begin{restatable}{claim}{claimCtoC}
    \label{CtoC}
    Suppose $x_1$ is the head of a clique $X \in C$ and let $Y \in C.$ Then, we may assume $\deg(x_1, Y) \leq r.$
  \end{restatable}

  With these claims at hand, we are ready to prove the main theorem. Recall that we assumed $G$ has no $r^{th}$ power of a Hamilton cycle with absolute discrepancy at least $\gamma n.$ We can assume $|f_R(\calT)| < \beta \ell,$ as otherwise we can find the desired $r^{th}$ power of a Hamilton cycle by Claim~\ref{tilingtoHr_claim}. Note that
  \begin{equation} \label{eq:sum_of_abcd}
    |A| + |B| + |C| + |D| = \ell / (r+1).
  \end{equation}
  Without loss of generality, we may assume
  \begin{equation}
    \label{bcad}
    |B| + |C| \geq |A| + |D|.
  \end{equation}

  First, we show that $A = \emptyset.$ Otherwise, consider some vertex $v$ of a clique in $A.$ By Claims \ref{AtoB} and \ref{AtoC}, $v$ can have at most $r-1$ edges toward any clique in $B \cup C$ and it can trivially have at most $r+1$ edges to any clique in $A \cup D.$ Using \eqref{eq:sum_of_abcd} and \eqref{bcad}, we get
  \[ \deg(v) \leq (r-1)(|B| + |C|) + (r+1)(|A| + |D|) \leq \frac{r}{r+1} \ell, \]
  which contradicts the degree assumption \eqref{degreeofR}. Hence, $A = \emptyset.$

  Note that 
  \[ f_R(\calT) = -{r+1 \choose 2}|B| + \left( -{r+1 \choose 2} + 2r\right) (|C| - |D|) = \left(- {r+1 \choose 2} + 2r\right) (|B| + |C| - |D|) - 2r |B| \le -2r |B|, \]
  where in the last inequality we used \eqref{bcad} and the fact that ${r+1 \choose 2} \ge 2r$ for $r \ge 3.$ By our assumption, it follows that $|B| < \beta \ell.$

  Next we show that $D = \emptyset.$ Indeed, suppose this is not the case and let $v$ be the head of some clique in $D.$ Then using Claim~\ref{DtoC}, \eqref{eq:sum_of_abcd} and \eqref{bcad}, we obtain
  \begin{align*}
    \deg(v) &\le (r+1)(|B| + |D|) + (r-1) |C| = \frac{r-1}{r+1} \ell + 2(|B| + |D|) \le \frac{r-1}{r+1} \ell + 2\left(\beta \ell + \frac{\ell}{2(r+1)}\right) \\ &\le \left( \frac{r}{r+1} + 2\beta\right) \ell < \delta(R),
  \end{align*}
  a contradiction. Thus, $D = \emptyset.$

  Observe that $|C| \ge \frac{\ell}{r+1} - |B| > 0.$ Now let $v$ be the head of a $(K_{r+1}, +)$-star in $C.$ By Claim~\ref{CtoC}, we can bound its degree as follows:
  \[ \deg(v) \le (r+1)|B| + r|C| = \frac{r}{r+1} \ell + |B| \le \left( \frac{r}{r+1} + \beta \right) \ell < \delta(R). \]
  Again, we have reached a contradiction so the proof is complete.
\end{proof}

\subsection*{Proofs of Claims~\ref{AtoB}--\ref{CtoC}}
In the following we will consider two $(r+1)$-cliques in $\calT$ which we denote by $X = \{x_1, \dots, x_{r+1}\}$ and $Y = \{y_1, \dots, y_{r+1}\}.$ Additionally, we assume that $x_1$ has at least $r$ edges towards vertices of $Y.$ Without loss of generality, we assume $x_1$ is connected to $y_1, y_2, \dots, y_r.$ From Claim \ref{common_neigh} it follows that there is a vertex $x'$ such that $X' = X \cup \{x'\}$ forms an $(r+2)$-clique. Similarly, there is a vertex $y'$ such that $Y' = Y \cup \{y'\}$ forms an $(r+2)$-clique and all the vertices $x_1, \dots x_{r+1}, x', y_1, \dots, y_{r+1}, y'$ are distinct. We construct two templates on $F = R[X' \stackrel{.}{\cup} Y']$ and show, for the cases mentioned in the claims, that these templates have different discrepancy. By Claim \ref{template_claim}, this contradicts the assumption that $G$ has no $r^{th}$ power of a Hamilton cycle with a large discrepancy. We start by defining four cycles which will be used in the templates.
\begin{align*}
  C_1 = (&x_2, x_3, \dots, x_{r+1}, x') \\
  C_2 = (&x_1, x_2, x_3, \dots, x_{r+1}, x') \\
  C_3 = (&y_1, y_2, y_3, \dots, y_{r+1}, y')\\
  C_4 = (&x_1, y_1, y_2, \dots, y_r, \\
        &y_{r+1}, y', y_1, y_2, \dots, y_{r-2}, y_{r-1}, \\
        &y_{r+1}, y', y_1, y_2, \dots, y_{r-2}, y_r, \\
        &\dots,\\
        &y_{r+1}, y', y_1, y_3, \dots, y_r, \\
        &y_{r+1}, y', y_2, y_3, \dots, y_r, \\
        &y_{r+1}, y', y_1, y_2, \dots, y_r)
\end{align*}
\\
Using these, we define two templates $\calF_1$ and $\calF_2$ as follows:
\begin{align*}
  \calF_1 &= \big((r+1) \times C_2, (r+1) \times C_3\big) \; \text{and} \\
  \calF_2 &= \big(C_1, r \times C_2, C_4\big),
\end{align*}
where we write $a \times C_i$ to indicate $a$ copies of $C_i.$ Note that $\calF_1$ contains each vertex in $F$ exactly $r+1$ times. Each vertex in $X' \setminus \{x_1\}$ appears once in $C_1$ and $C_2,$ so $\calF_2$ contains each of these vertices $r+1$ times. $C_4$ contains $x_1$ once and each of the vertices in $Y'$ exactly $r+1$ times which gives 
\begin{equation}
  \label{c4_len}
  |C_4| = 1 + (r+2)(r+1) = r^2 + 3r + 3,
\end{equation}
which will be useful later on. Additionally, it is easy to see that $\calF_2$ contains each vertex in $F$ exactly $r+1$ times. Therefore, if $\calF_1$ and $\calF_2$ have different discrepancies, we reach a contradiction by Claim \ref{template_claim}. In the following claims we show this is true for several cases of interest. When calculating the discrepancy of a particular $r^{th}$ power of a cycle, we will mostly use the following recipe. As the $(r+1)$-cliques we consider are highly structured, most of the edge values under consideration are known given the types of cliques of $X$ and $Y.$ More precisely, we find a small subset of edges $E'$ such that all edges in $E(C_i^r) \setminus E'$ have the same color $c$. Then, we can calculate the discrepancy of $C_i^r$ as
\[f_R(C_i^r) = c \big(r|C_i| - |E'| \big) + f_R(E').\]
Additionally, observe that
\begin{equation}
  \label{diff}
  f_R(\calF_1) - f_R(\calF_2) = -f_R(C_1^r) + f_R(C_2^r) + (r+1) f_R(C_3^r) - f_R(C_4^r).
\end{equation}

Finally, we proceed to prove the individual claims.
\claimAtoB*
\begin{proof}
  Suppose $\deg(x_1, Y) \geq r$ and let $\calF_1, \calF_2$ be defined as above. Applying Claim \ref{cliques_cor} to the clique $X',$ we get that all edges from $x'$ to $X$ have the same color $f_R(x', x_1).$ Analogously, all edges from $y'$ to $Y$ have color $f_R(y', y_1)$ and all edges from $x_1$ to $Y \setminus \{y_{r+1}\}$ have color $f_R(x_1, y_1).$ We calculate the discrepancies of $C_i^r, 1 \le i \le 4$ following our recipe described above. For example, $f(C_1^r)$ is calculated as follows. Observe that $C_1^r$ contains $r+1$ vertices and every vertex has degree $2r$ in $C_1^r$ so there are $r(r+1)$ edges. All edges not incident to $x'$ are between two vertices of $X$ and thus have color $+1$ because $X$ forms a copy of $K_{r+1}^+.$ There are $r(r+1) - 2r$ such edges and the remaining $2r$ edges have the same color $f_R(x', x_1)$ as shown. Therefore, we obtain
  \[ f_R(C_1^r) = ((r+1)r - 2r) \cdot 1 + 2r \cdot f_R(x', x_1) = r^2 - r + 2r f_R(x', x_1). \]
  A similar analysis yields the following:
  \begin{align*}
    f_R(C_2^r) &= (r+2)r - 2r + 2r f_R(x', x_1) = r^2 + 2r f_R(x', x_1) \\
    f_R(C_3^r) &= -((r+2)r - 2r) + 2rf_R(y', y_1) = -r^2 + 2r f_R(y', y_1) \\
    f_R(C_4^r) &= -\big(r|C_4| - 2r - 2r(r+1)\big) + 2r f_R(x_1, y_1) + 2r(r+1) f_R(y', y_1) \\
          &= -r^3 - r^2 + r + 2rf_R(x_1, y_1) + 2r(r+1) f_R(y', y_1).
  \end{align*}
  Plugging these values into \eqref{diff}, we obtain:
  \begin{align*}
    f_R(\calF_1) - f_R(\calF_2) = -2r f_R(x_1, y_1) \neq 0.
  \end{align*}
	We are done by Claim \ref{template_claim}.
\end{proof}

\claimAtoC*
\begin{proof}
  Assume $\deg(x_1, Y) \geq r$ and define $\calF_1, \calF_2$ as above. By Claim \ref{cliques_cor}, all edges from $x'$ to $X$ have the same color $f_R(x', x_1).$ By Claim \ref{cliques_cor}, $Y'$ is a $(K_{r+2}, +)$-star with its head in $Y.$ Note that the values of $f_R(C_1^r)$ and $f_R(C_2^r)$ are as in the previous claim; we also calculate $f_R(C_3^r)$:
  \begin{align*}
    f_R(C_1^r) &= r^2 - r + 2r f_R(x', x_1) \\
    f_R(C_2^r) &= r^2 + 2r f_R(x', x_1) \\
    f_R(C_3^r) &= -\big((r+2)r - 2r\big) + 2r = -r(r-2).
  \end{align*}
  
  Now, we consider two cases:
  \begin{enumerate}[label=(\alph*)]
    \item $y_{r+1}$ is the head of $Y.$ \\
    Applying Claim \ref{cliques_cor} to $(Y \setminus \{y_{r+1}\}) \cup \{x_1\},$  we conclude that all edges from $x_1$ to $y_i, \; i \in [r]$ are of color $f_R(x_1, y_1).$ Using this, we obtain:
    \begin{align*}
      f_R(C_4^r) &= -(r|C_4| - 2r - 2r(r+1)) + 2rf_R(x_1, y_1) + 2r(r+1)\\
            &= -r^3 + r^2 + 3r + 2r f_R(x_1, y_1).
    \end{align*}
    Substituting into \eqref{diff}, we have
    \[ f_R(\calF_1) - f_R(\calF_2) = -2r f_R(x_1, y_1) \neq 0. \] 

    \item $y_{r+1}$ is not the head of $Y.$ \\
    Applying Claim \ref{cliques_cor} to $(Y \setminus \{y_{r+1}\}) \cup \{x_1\},$ we get that all edges in $C_4^r$ incident to the head of $Y$ have value $+1,$ while all other edges have value $-1.$ From this we have:
    \[ f_R(C_4^r) = -(r|C_4| - 2r(r+1)) + 2r(r+1) = -r^3 + r^2 + r \]
    and
    \[ f_R(\calF_1) - f_R(\calF_2) = 2r \neq 0\]
  \end{enumerate}
  In both cases, the proof is finished by Claim \ref{template_claim}. 
\end{proof}

\claimDtoC*
\begin{proof}
  Again, suppose $\deg(x_1, Y) \geq r$ and define $\calF_1$ and $\calF_2$ as above. By Claim \ref{cliques_cor}, $Y'$ is a $(K_{r+2}, +)$-star with its head in $Y.$ We have
  \[ f_R(C_1^r) = r(r+1),
    f_R(C_2^r) = r(r-2) \text { and }
    f_R(C_3^r) = -r(r-2). \]

  Note that the same edges of $C_4^r$ are known as those in Claim \ref{AtoC}. Again, we consider two cases:
  \begin{enumerate}[label=(\alph*)]
    \item $y_{r+1}$ is the head of $Y.$ \\
    Applying Claim \ref{cliques_cor} to $(Y \setminus \{y_{r+1}\}) \cup \{x_1\},$ we conclude that all edges from $x_1$ to $y_i, \; i \in [r]$ are of color $f_R(x_1, y_1).$ So, we get 
    \[ f_R(C_4^r) = -r^3 + r^2 + 3r + 2r f_R(x_1, y_1). \]
    Substituting into \eqref{diff}, we have
    \[ f_R(\calF_1) - f_R(\calF_2) = -4r - 2r f_R(x_1, y_1) \neq 0. \] 

    \item $y_{r+1}$ is not the head of $Y.$ \\
    We apply Claim \ref{cliques_cor} to $(Y \setminus \{y_{r+1}\}) \cup \{x_1\}$ and obtain that all edges in $C_4^r$ incident to the head of $Y$ have value $1,$ while all other edges have value $-1.$ From this we have:
    \[ f_R(C_4^r) = -r^3 + r^2 + r \]
    and
    \[ f_R(\calF_1) - f_R(\calF_2) = -2r \neq 0\]
  \end{enumerate}
  Again, we are done by Claim \ref{template_claim}.
\end{proof}

\claimCtoC*
\begin{proof}
  Let $x_1$ be the head of a clique in $C,$ let $Y$ be a clique in $C$ and assume that $\deg_R(x_1, Y) = r+1.$ Clearly, $Y$ is not the clique containing $x_1.$ Since $x_1$ is connected to all the vertices in $Y$ we may assume that $y_1$ is the head of $Y$ and use the same template as before. Similarly as in the proof of Claim~\ref{DtoC}, we have that $X'$ and $Y'$ are $(K_{r+2}, +)$-stars with heads $x_1$ and $y_1,$ respectively. Since $\{x_1\} \cup Y$ forms a clique and $Y$ is a $(K_{r+1},+)$-star, it follows that $\{x_1\} \cup Y$ is a $(K_{r+2}, +)$-star with head $y_1.$ This implies that $f_R(x_1, y_1) = 1$ and $f_R(x_1, y_i) = -1, 2 \le i \le r+1.$ As in the proof of Claim~\ref{DtoC}, but with opposite signs, we have:
  \[ f_R(C_1^r) = -r(r+1), \, f_R(C_2^r) = -r(r-2) \text{ and } f_R(C_3^r) = -r(r-2). \]
  Note that all edges of $C_4^r$ have value $-1$ apart from all the edges incident to $y_1.$ Therefore,
  \[ f_R(C_4^r) =  -r(|C_4^r| - 2(r+1)) + 2r(r+1) = -r^3 + r^2 + r. \]
  Substituting into \eqref{diff}, we get
  \[ f_R(\calF_1) - f_R(\calF_2) = r(r+1) - r(r-2) - (r+1)r(r-2) + r^3 - r^2 - r = 4r \]
  and we are done by Claim~\ref{template_claim}.
\end{proof}

\subsection*{Acknowledgments}
The author would like to thank Milo\v{s} Truji\'{c} for introducing him to the problem, helpful discussions and
many useful comments. Additionally, the author is grateful for the comments of the anonymous referees
which greatly improved the presentation of this paper, especially the referee who noticed a mistake in the
original argument for the case $r = 3.$

%%%%%%%%%%%%%%%%%%%%%%%%%%%%%%
{\small
\bibliographystyle{plain}
\bibliography{references}
}

\appendix
% !TeX root = main.tex
% !TEX root = main.tex

\section{Appendix}
Here we present the missing proofs of Claims~\ref{cl:conn_paths}~and~\ref{cl:adding_vertices_to_paths}. We begin with an auxiliary claim which will be used in the proofs.

\begin{claim} \label{cl:finding_vertices}
Let $C_{-r+1}, C_{-r+2}, \dots, C_{t+r}$ be a sequnce of clusters forming the $r^{th}$ power of a path in $R$ with $t = O(r^3).$ Furthermore, let $U_i \subseteq C_i,$ for $-r+1 \le i \le t+r,$ be given subsets of size at least $(d-\varepsilon)^{3r} m / 8.$ Then, there exist vertices $p_1, \dots, p_t$ with $p_i \in U_i, \forall i \in [t],$ forming the $r^{th}$ power of a path in $G$ such that
  \begin{equation} \label{eq:extensible_endpoints}
    \begin{aligned}
      |N(\{p_1, p_2, \dots, p_j\}) \cap U_{-r+j}| &\ge (d-\varepsilon)^{2r} |U_{-r+j}|/2, \forall j \in [r] \text{ and }\\
      |N(\{p_t, p_{t-1}, \dots, p_{t+1-j}\}) \cap U_{t+r+1-j}| &\ge (d - \varepsilon)^{2r} |U_{t+r+1-j}|/2, \forall j \in [r] .
    \end{aligned}
  \end{equation}
\end{claim}
\begin{proof}
  We will choose the vertices $p_1, \dots, p_t$ one by one. For this purpose, we maintain sets $H_{i,j}$ from which the vertices will be selected starting with $H_{0,j} = U_j, \, -r+1 \le j \le r+t.$ Assume we have selected vertices $p_1, \dots, p_{i-1}$ and that $|H_{i-1, i}| > 2r \varepsilon |C_i|,$ which will follow from the construction. Then we select a vertex $p_i$ from $H_{i-1, i}$ satisfying:
  \[ \deg(p_i, H_{i-1, j}) > (d - \varepsilon) |H_{i-1, j}| \text{ for all } j \neq i, |j-i| \le r. \]
  Provided $|H_{i-1, j}| > \varepsilon |C_j|,$ by $(\varepsilon, d)$-regularity, this holds for all but at most $2r \varepsilon |C_i| < |H_{i-1,i}|$ vertices in $H_{i-1,i},$ so we can choose such a vertex. We then update the sets as follows:
  \[ H_{i, j} = \begin{cases}
    H_{i-1, j} \cap N(p_i), &\text{if } 1 \le |j-i| \le r,\\
    H_{i-1, j} \setminus \{p_i\}, &\text{otherwise}. 
  \end{cases} \]
  We need to argue that the sets $H_{i, j}$ are large enough for the above arguments to hold. This follows from the fact that for each $j, \, 1 \le j \le t,$ the set $H_{0, j}$ shrinks by a factor of $d-\varepsilon$ at most $2r$ times and we remove from it single vertices $O(r^3)$ times. Since we can assume $(d-\varepsilon)^5 > 32r \varepsilon$ and $m$ is large compared to $r,$ we obtain that $|H_{i, j}| \ge (d-\varepsilon)^{2r} |U_j| / 2 > 2r\varepsilon |C_j|,$ for all $i, j.$

  By construction, the vertices $p_1, \dots, p_t$ form the $r^{th}$ power of a path while the property~\eqref{eq:extensible_endpoints} follows from $|H_{t, j}| > (d-\varepsilon)^{2r} |U_j| / 2$ for $j < 0$ and $j > t,$ thus finishing the proof.
\end{proof}

\claimconnectingpaths*
\begin{proof}
  We construct these paths one by one in order. We show how to construct a path $\calP_1$ satisfying the desired properties. The remaining paths can be constructed analogously, the only difference being that we cannot use the vertices in the previously constructed paths. This can be easily guaranteed as at any point the number of used vertices is only a constant compared to $m,$ the size of each cluster.

  First, we determine the sequence of clusters from which we choose the vertices of $\calP_1.$ This sequence of clusters will form the $r^{th}$ power of a (not necessarily simple) path in $R.$ To achieve this, we find a sequence $L_1, L_2, \dots, L_t$ of $(r+1)$-cliques in $R$ such that:
  \begin{enumerate}[label=\alph*)]
    \item \label{seq1} $L_1 = K_1, \, L_t = K_2,$
    \item $|L_{i+1} \cap L_i| = r,$ for all $1 \le i \le t-1,$
    \item \label{seq3} $t = O(r^2).$
  \end{enumerate}
  For this purpose, for any two cliques $L, L'$ in $R$ and any cluster $C \in V(R),$ we define its weight with respect to $L, L'$ as $w_{(L, L')}(C) = \deg_R(C, L) + \deg_R(C, L').$ We will need the following simple claim.
  \begin{claim} \label{cl:labels}
    For any two $(r+1)$-cliques $L, L'$ in $R$ there are at least $\frac{\eta'}{r+1} \ell$ clusters $C \in V(R)$ such that \linebreak $w_{(L, L')}(C) \ge 2r+1.$
  \end{claim}
  \begin{proof}
    By \eqref{eq:degR}, we have
    \[ \sum_{C \in V(R)} w_{(L, L')}(C) \ge 2(r+1) \delta(R) \ge 2(r + \eta') \ell. \]
    The contribution to the above sum of the clusters $C$ with $w_{(L,L')}(C) \le 2r$ is at most $2r \ell.$ Since for any $C,$ $w_{(L, L')}(C) \le 2r+2,$ the statement follows.
  \end{proof}
  We construct the sequence $L_1, \dots, L_t$ of $(r+1)$-cliques in two phases. First, we construct two sequences of $(r+1)$-cliques $A_1, A_2, \dots, A_{t_1}$ and $B_1, B_2, \dots, B_{t_2}$ with the following properties:
  \begin{enumerate}[label=\arabic*)]
    \item $A_1 = K_1, B_1 = K_2,$
    \item $|A_i \cap A_{i+1}| = r, |B_j \cap B_{j+1}| = r,$ for all $1 \le i \le t_1-1, 1 \le j \le t_2-1,$
    \item \label{done_sequences} $\deg_R(C, A_{t_1}) \ge r \text{ for all } C \in B_{t_2},$ or $\deg_R(C, B_{t_2}) \ge r \text{ for all } C \in A_{t_1}.$
  \end{enumerate}
  These sequences can be construted as follows. Set $A_1 = K_1, B_1 = K_2$ and suppose we have so far constructed $A_1, \dots, A_i$ and $B_1, \dots, B_j$ and \ref{done_sequences} does not hold. By Claim~\ref{cl:labels}, there exists a cluster $C \not\in A_i \cup B_j$ such that $w_{(A_i, B_j)}(C) \ge 2r+1.$ Hence $d_R(C, A_i) \ge r+1$ and $d_R(C, B_j) \ge r$ or vice-versa. Assume the former holds, the other case being handled analogously. As \ref{done_sequences} is not satisfied, there is a cluster $C' \in A_i$ with $\deg_R(C', B_j) \le r-1.$ To get $A_{i+1},$ from $A_i$ we remove $C'$ and add $C.$ Note that we maintain the property that the last cliques in the two sequences, $A_{i+1}$ and $B_j,$ are disjoint and in each step we increased the number of edges between them so in $O(r^2)$ steps \ref{done_sequences} is satisfied.

  Now we construct a sequence of $(r+1)$-cliques $D_0, \dots, D_{r+1}$ forming a transition from $A_{t_1}$ to $B_{t_2}.$ Assume that $\deg_R(C, A_{t_1}) \ge r$ for all $C \in B_{t_2},$ the other case being analogous.
  We will construct these cliques such that for all $0 \le i \le r+1,$ $|D_i \cap A_{t_1}| = r+1-i$ and $|D_i \cap B_{t_2}| = i.$ Set $D_0 = A_{t_1},$ let $0 \le i \le r$ and suppose we have already constructed $D_0, \dots, D_i.$ Let $C$ be an arbitrary cluster in $B_{t_2} \setminus D_i.$ By assumption, there is at most one cluster in $D_i \setminus B_{t_2} \subseteq A_{t_1}$ not adjacent to $C.$ We remove this cluster, or if none exists, an arbitrary cluster in $D_i \setminus B_{t_2}$ from $D_i$ and add $C$ to obtain $D_{i+1}.$

By construction, the clique sequence
\[ A_1, A_2, \dots, A_{t_1}, D_1, D_2, \dots, D_r, B_{t_2}, B_{t_2-1}, \dots, B_1 \]
satisfies \ref{seq1}--\ref{seq3}.
We denote this sequence by $L_1, \dots L_t,$ where $t = t_1 + t_2 + r = O(r^2).$
To get the sequence of clusters from which we will choose the vertices of $\calP_1,$ we do the following. Recall that $L_1 = V^1.$ We start with $V^1_1, V^1_2, \dots, V^1_{r+1}.$ Then we start another cycle $V^1_1, V^1_2, \dots$ until we reach the cluster before the unique cluster in $L_1 \setminus L_2.$ The next cluster in the sequence is then the unique cluster in $L_2 \setminus L_1.$ To ensure the sequence of clusters forms the $r^{th}$ power of a path we make a cycle through the clusters in $L_2$ and part of another cycle before reaching the cluster in $L_2 \setminus L_3.$ We append the unique cluster in $L_3 \setminus L_2$ and so on. We proceed in a similar fashion until we reach the point where the last $r+1$ clusters form the clique $K_2.$ Therefore, we have a seqence of clusters $C_1, \dots, C_{t'}$ with $t' = O(r^3)$ such that $C_i = V^1_i,$ for $i \in [r+1],$ and the last $r+1$ clusters are the clusters of $K_2$ in some order.

However, we want the last $r+1$ clusters to be $V^2_1, \dots, V^2_{r+1}$ in precisely this order. To achieve this, it is enough to show that if $U_1, \dots, U_{r+1}$ is an arbitrary permutation of $V^2_1, \dots, V^2_{r+1}$ and $1 \le i < j \le r+1,$ then we can create a sequence of clusters forming the $r^{th}$ power of a path in $R$ which starts with $U_1, \dots, U_{r+1}$ and ends with $U_1, \dots, U_{i-1}, U_j, U_{i+1}, \dots, U_{j-1}, U_i, U_{j+1}, \dots, U_{r+1}.$ By \eqref{eq:common_neigh}, there is a cluster $C \in N_R(\{U_1, \dots, U_{r+1}\}).$ The desired sequence is then given as follows:
\begin{align*}
  \begin{gathered}
  U_1, \dots, U_{r+1},\\
  U_1, \dots, U_{i-1}, C, U_{i+1}, \dots, U_{r+1},\\
  U_1, \dots, U_{i-1}, C, U_{i+1}, \dots, U_{j-1}, U_i, U_{j+1}, \dots, U_{r+1},\\
  U_1, \dots, U_{i-1}, U_j, U_{i+1}, \dots, U_{j-1}, U_i, U_{j+1}, \dots, U_{r+1}.
  \end{gathered}
\end{align*}

As we can get the desired permutation using at most $r$ of these swaps, we obtain a sequence of clusters $C_1, \dots, C_{t''}$ which form the $r^{th}$ power of a path in $R$ of length $t'' = O(r^3)$ such that for all $j \in [r+1],$ $C_j = V^1_j$ and $C_{t''-r-1+j} = V^2_j.$ We further denote
\[ C_0 = V^1_{r+1}, C_{-1} = V^1_r, \dots, C_{-r+1} = V^1_2, \; \text{ and } \; C_{t''+1} = V^2_1, C_{t''+2} = V^2_2, \dots, C_{t''+r} = V^2_r. \]
Now we can find $\calP_1$ by applying Claim~\ref{cl:finding_vertices} with $U_j = C_j, \forall j, -r+1 \le j \le t''+r.$

When constructing the path $\calP_k,$ the only difference to the above argument is that we remove the set of vertices used in the previous paths from the sets $U_i.$ Since we remove at most $\ell \cdot O(r^3) < m/2$ vertices, Claim~\ref{cl:finding_vertices} is still applicable.
\end{proof}

\claimaddingvertices*
\begin{proof}
  For convenience of notation we only show how to extend $\calP_1$ and the other paths are extended analogously. Recall that $|A_1| \le \varepsilon_1 m' = \frac{100 r^3 \varepsilon}{\eta} m'.$ We embed the vertices in $A_1$ one by one while maintaining property \ref{enum:extensible_endpoints}. Slightly abusing notation, let $p_{-r+1}, p_{-r+2}, \dots, p_0$ denote the last $r$ vertices of $\calP_1,$ where by construction $p_i \in V^2_{r+1+i}, -r+1 \le i \le 0,$ and suppose we wish to extend $\calP_1$ to additionally contain a vertex $v \in A_1$ which, by definition, satisfies
  \begin{equation} \label{eq:v_is_assigned}
    \deg(v, V_j^2) \ge (d - \varepsilon)|V_j^2|, \, \forall j \in [r+1].
  \end{equation}
  We aim to extend $\calP_1$ by going around the clusters in $K_2$ twice, then embedding $v$ and then going around the clusters in $K_2$ once again. To achieve this, we will apply Claim~\ref{cl:finding_vertices}. We set $C_{(r+1)k + j} = V^2_j$ for $-1 \le k \le 3$ and $1 \le j \le r+1.$ Additionally we set, 
  \[ U_{(r+1)k + j} = \begin{cases}
    \big(C_{(r+1)k+j} \cap N(\{ p_i \, \vert \, -r+1 \le i \le 0, p_i \not\in C_{(r+1)k+j}\}\big) \setminus \mathrm{Used}, &\text{for } k \in \{-1, 0\}, j \in [r+1],\\
    \big(C_{(r+1)k+j} \cap N(v)\big) \setminus \mathrm{Used}, &\text{for } k \in \{1, 2\}, j \in [r+1],\\
    \big(C_{(r+1)k+j} \setminus \mathrm{Used}, &\text{for } k = 3, j \in [r+1].
    \end{cases}
  \]
  where $\mathrm{Used}$ denotes the set of vertices already used in all of the paths so far plus the new vertex $v.$ Note that $|\mathrm {Used} \cap C_{(r+1)k+j}| \le \ell \cdot O(r^3) + 3\varepsilon m' < m/8,$ for all $-1 \le k \le 3, 1 \le j \le r+1.$ Using this, the fact that $\calP_1$ satisfies \ref{enum:extensible_endpoints} and \eqref{eq:v_is_assigned}, it follows that
  \[ 
    |U_{(r+1)k + j}| \ge \begin{cases}
      (d-\varepsilon)^{2r} m / 8 &\text{for } -1 \le k \le 2, 1 \le j \le r+1, \\
      m/2 &\text{for } k = 3, 1 \le j \le r+1.
    \end{cases}
  \]

  Finally, let $q_1, \dots, q_{3(r+1)}$ be the vertices obtained by applying Claim~\ref{cl:finding_vertices} for the sets $C_i, U_i, \, -r+1 \le i \le 4r+3.$ We extend the path $\calP_i$ with the vertices $q_1, q_2, \dots, q_{2(r+1)}, v, q_{2(r+1)+1}, \dots, q_{3(r+1)}$ in this order. Using $|U_{3(r+1)+j}| \ge m/2$ for $1 \le j \le r+1$ and \eqref{eq:extensible_endpoints}, we see that the new new path still satisfies \ref{enum:extensible_endpoints}. By definition of $U_{k(r+1)+j}$ for $1 \le k \le 2,$ it follows that $v q_{k(r+1)+j} \in E(G)$ for $1 \le k \le 2, j \in [r+1].$ Since $q_1, \dots, q_{3(r+1)}$ form the $r^{th}$ power of a path, it follows that $\calP_1$ is still the $r^{th}$ power of a path, completing the proof.

\end{proof}

\end{document}